\documentclass[letter,12pt]{article}
\usepackage{skak}
\usepackage{xkeyval}
\usepackage{chessfss}
\usepackage{xifthen}
\usepackage[latin1]{inputenc}
\usepackage[english]{babel}
\usepackage{latexsym}
\usepackage{amsmath} 
\usepackage{fullpage}
\usepackage{textcomp}
\usepackage{graphicx}
\usepackage{amssymb}
\usepackage{float}
\usepackage{epstopdf}					
\usepackage[pdftex]{hyperref}					
 \usepackage{pdflscape}
 \usepackage{mathrsfs}
\usepackage{tikz}
 \usepackage{pifont}
 \usepackage{multicol}
\definecolor{1}{rgb}{1,0.2,0.3}
\definecolor{2}{rgb}{0.1,0.3,0.5}
\definecolor{3}{rgb}{1,1,0}
\definecolor{4}{rgb}{0.2,0.6,0}

\newtheorem{eg}{Example}

\newcommand{\N}{\mathbb{N}}

\usepackage{fancyhdr}
\begin{document}
\begin{center}{\textbf{ \begin{large}The Mutando of Insanity \end{large} \begin{small} \\ by \'Erika. B. Rold\'an Roa  \end{small} } \\  \line(1,0){530}}
\end{center}
Puzzles based on coloured cubes and other coloured geometrical figures have a long history in the recreational mathematical literature. Martin Gardner wrote about them in: Chapter 16 of \textit{New Mathematical Diversions}, Chapter 16 of \textit{Mathematical Magic Show }, and Chapter 6 of \textit{Fractal Music, Hypercards and More Mathematical Recreations from Scientific American Magazine} (\cite{gardner1995new, gardner1977mathematical, gardner1992fractal}). One of the most commercially famous of these puzzles is the Instant Insanity that consists of four cubes. Their faces are coloured with four different colours in such a way that each colour is present in each one of the four cubes. To solve the puzzle, one needs to stack the cubes in a tower in such a way that each one of the colours appears exactly once in the four long faces of the tower. The main purpose of this paper is to study the combinatorial richness of a mathematical model of this puzzle by analysing all possible ways of colouring cubes to form a puzzle analogous to the Instant Insanity. We have done this analysis for $n$ cubes and $n$ colours for $n=4, 5, 6$. This combinatorial analysis allowed us to design the Mutando of Insanity, a puzzle that we presented in a talk at the G4G 12.

\section*{Introduction}  
In the book \textit{A Lifetime of Puzzles}, Rik van Grol \cite{demaine2008lifetime} wrote a paper about the Insanity kind of puzzles. In this paper, he gave a wonderful account of the history of the "Instant Insanity" puzzle and related puzzles. The majority of the Insanity puzzles consist of $n$ cubes, with a fixed $n \in \{4,5,6\}$, and the six faces of each cube coloured with one of $n$ colours in such a way that each one of the four colours appears in each cube. To solve the puzzle one has to find a way of stacking the $n$ cubes in a tower (a $n \times 1 \times 1$ prism) in which the $n$ colours appear in each one of the long faces of the tower. As an example, for $n=4$ the way in which the cubes of The Instant Insanity puzzle are coloured is depicted in the next figure (they are not the original colours of the Instant Insanity but the structure is exactly the same):
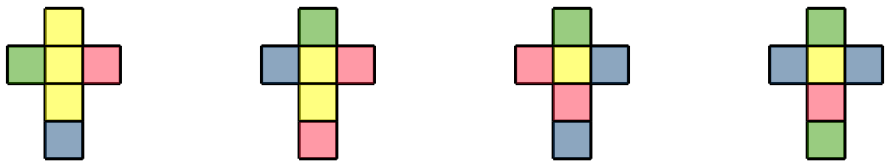
\begin{center}
 \begin{multicols}{4}

\begin{tikzpicture}[thick,scale=.5, ]
   \coordinate (A1) at (0,0);
   \coordinate (A2) at (1, 0);
    \coordinate (A3) at (0, 1);
    \coordinate (A4) at (1, 1);
    \coordinate (B1) at (0, 2);
    \coordinate (B2) at (1, 2);
    \coordinate (B3) at (0, 3);
    \coordinate (B4) at (1, 3);
   \coordinate (C1) at (0, 4);
    \coordinate (C2) at (1, 4);
   \coordinate(D1) at (-1,2);
   \coordinate(D2) at (-1,3);
  \coordinate(D3) at (2,2);
 \coordinate(D4) at (2,3);

   \draw[very thick] (A1) -- (A2);
    \draw[very thick] (A1) -- (A3);
    \draw[very thick] (A2) -- (A4);
    \draw[very thick] (A3) -- (A4);
   \draw[very thick] (A3) -- (B1);
   \draw[very thick] (A4) -- (B2);
   \draw[very thick] (B1) -- (B2);
   \draw[very thick] (B3) -- (C1);
   \draw[very thick] (B4) -- (C2);
   \draw[very thick] (B3) -- (B1);
   \draw[very thick] (B4) -- (B2);
   \draw[very thick] (B3) -- (B4);
 \draw[very thick] (C1) -- (C2);    
  \draw[very thick] (D1) -- (D2); 
 \draw[very thick] (D3) -- (D4);
  \draw[very thick] (B1) -- (D1);
  \draw[very thick] (B2) -- (D3);
  \draw[very thick] (B3) -- (D2);
  \draw[very thick] (B4) -- (D4);

\draw[fill=2,opacity=0.5]  (A1) -- (A2) -- (A4) -- (A3);
\draw[fill=3,opacity=0.5]  (A3) -- (A4) -- (B2) -- (B1);  
\draw[fill=3,opacity=0.5]  (B1) -- (B2) -- (B4) -- (B3);
\draw[fill=3,opacity=0.5]  (B3) -- (B4) -- (C2) -- (C1);
\draw[fill=4,opacity=0.5]  (D1) -- (D2) -- (B3) -- (B1);
\draw[fill=1 ,opacity=0.5]  (D3) -- (D4) -- (B4) -- (B2);
\end{tikzpicture}

\begin{tikzpicture}[thick,scale=.5,]
   \coordinate (A1) at (0,0);
   \coordinate (A2) at (1, 0);
    \coordinate (A3) at (0, 1);
    \coordinate (A4) at (1, 1);
    \coordinate (B1) at (0, 2);
    \coordinate (B2) at (1, 2);
    \coordinate (B3) at (0, 3);
    \coordinate (B4) at (1, 3);
   \coordinate (C1) at (0, 4);
    \coordinate (C2) at (1, 4);
   \coordinate(D1) at (-1,2);
   \coordinate(D2) at (-1,3);
  \coordinate(D3) at (2,2);
 \coordinate(D4) at (2,3);

   \draw[very thick] (A1) -- (A2);
    \draw[very thick] (A1) -- (A3);
    \draw[very thick] (A2) -- (A4);
    \draw[very thick] (A3) -- (A4);
   \draw[very thick] (A3) -- (B1);
   \draw[very thick] (A4) -- (B2);
   \draw[very thick] (B1) -- (B2);
   \draw[very thick] (B3) -- (C1);
   \draw[very thick] (B4) -- (C2);
   \draw[very thick] (B3) -- (B1);
   \draw[very thick] (B4) -- (B2);
   \draw[very thick] (B3) -- (B4);
 \draw[very thick] (C1) -- (C2);    
  \draw[very thick] (D1) -- (D2); 
 \draw[very thick] (D3) -- (D4);
  \draw[very thick] (B1) -- (D1);
  \draw[very thick] (B2) -- (D3);
  \draw[very thick] (B3) -- (D2);
  \draw[very thick] (B4) -- (D4);

\draw[fill=1,opacity=0.5]  (A1) -- (A2) -- (A4) -- (A3);
\draw[fill=3,opacity=0.5]  (A3) -- (A4) -- (B2) -- (B1);  
\draw[fill=3,opacity=0.5]  (B1) -- (B2) -- (B4) -- (B3);
\draw[fill=4,opacity=0.5]  (B3) -- (B4) -- (C2) -- (C1);
\draw[fill=2,opacity=0.5]  (D1) -- (D2) -- (B3) -- (B1);
\draw[fill=1,opacity=0.5]  (D3) -- (D4) -- (B4) -- (B2);

\end{tikzpicture}

\begin{tikzpicture}[thick,scale=.5,]
   \coordinate (A1) at (0,0);
   \coordinate (A2) at (1, 0);
    \coordinate (A3) at (0, 1);
    \coordinate (A4) at (1, 1);
    \coordinate (B1) at (0, 2);
    \coordinate (B2) at (1, 2);
    \coordinate (B3) at (0, 3);
    \coordinate (B4) at (1, 3);
   \coordinate (C1) at (0, 4);
    \coordinate (C2) at (1, 4);
   \coordinate(D1) at (-1,2);
   \coordinate(D2) at (-1,3);
  \coordinate(D3) at (2,2);
 \coordinate(D4) at (2,3);

   \draw[very thick] (A1) -- (A2) ;
    \draw[very thick] (A1) -- (A3) ;
    \draw[very thick] (A2) -- (A4) ;
    \draw[very thick] (A3) -- (A4);
   \draw[very thick] (A3) -- (B1);
   \draw[very thick] (A4) -- (B2);
   \draw[very thick] (B1) -- (B2);
   \draw[very thick] (B3) -- (C1);
   \draw[very thick] (B4) -- (C2);
   \draw[very thick] (B3) -- (B1);
   \draw[very thick] (B4) -- (B2);
   \draw[very thick] (B3) -- (B4);
 \draw[very thick] (C1) -- (C2);    
  \draw[very thick] (D1) -- (D2); 
 \draw[very thick] (D3) -- (D4);
  \draw[very thick] (B1) -- (D1);
  \draw[very thick] (B2) -- (D3);
  \draw[very thick] (B3) -- (D2);
  \draw[very thick] (B4) -- (D4);

\draw[fill=2,opacity=0.5]  (A1) -- (A2) -- (A4) -- (A3);
\draw[fill=1,opacity=0.5]  (A3) -- (A4) -- (B2) -- (B1);  
\draw[fill=3,opacity=0.5]  (B1) -- (B2) -- (B4) -- (B3);
\draw[fill=4,opacity=0.5]  (B3) -- (B4) -- (C2) -- (C1);
\draw[fill=1,opacity=0.5]  (D1) -- (D2) -- (B3) -- (B1);
\draw[fill=2,opacity=0.5]  (D3) -- (D4) -- (B4) -- (B2);
\end{tikzpicture}

\begin{tikzpicture}[thick,scale=.5, ]
   \coordinate (A1) at (0,0);
   \coordinate (A2) at (1, 0);
    \coordinate (A3) at (0, 1);
    \coordinate (A4) at (1, 1);
    \coordinate (B1) at (0, 2);
    \coordinate (B2) at (1, 2);
    \coordinate (B3) at (0, 3);
    \coordinate (B4) at (1, 3);
   \coordinate (C1) at (0, 4);
    \coordinate (C2) at (1, 4);
   \coordinate(D1) at (-1,2);
   \coordinate(D2) at (-1,3);
  \coordinate(D3) at (2,2);
 \coordinate(D4) at (2,3);

   \draw[very thick] (A1) -- (A2);
    \draw[very thick] (A1) -- (A3);
    \draw[very thick] (A2) -- (A4);
    \draw[very thick] (A3) -- (A4);
   \draw[very thick] (A3) -- (B1);
   \draw[very thick] (A4) -- (B2);
   \draw[very thick] (B1) -- (B2);
   \draw[very thick] (B3) -- (C1);
   \draw[very thick] (B4) -- (C2);
   \draw[very thick] (B3) -- (B1);
   \draw[very thick] (B4) -- (B2);
   \draw[very thick] (B3) -- (B4);
 \draw[very thick] (C1) -- (C2);    
  \draw[very thick] (D1) -- (D2); 
 \draw[very thick] (D3) -- (D4);
  \draw[very thick] (B1) -- (D1);
  \draw[very thick] (B2) -- (D3);
  \draw[very thick] (B3) -- (D2);
  \draw[very thick] (B4) -- (D4);

\draw[fill=4,opacity=0.5]  (A1) -- (A2) -- (A4) -- (A3);
\draw[fill=1,opacity=0.5]  (A3) -- (A4) -- (B2) -- (B1);  
\draw[fill=3,opacity=0.5]  (B1) -- (B2) -- (B4) -- (B3);
\draw[fill=4,opacity=0.5]  (B3) -- (B4) -- (C2) -- (C1);
\draw[fill=2,opacity=0.5]  (D1) -- (D2) -- (B3) -- (B1);
\draw[fill=2,opacity=0.5]  (D3) -- (D4) -- (B4) -- (B2);
\end{tikzpicture}
\end{multicols}
\end{center}

 For the rest of the paper, when we refer to an Insanity puzzle we mean a collection of $n$ coloured cubes with $n$ colours as we described it in the last paragraph. For each fixed $n \in \{4,5,6\}$, there are a lot of ways of colouring one cube with $n$ colours and even more ways to select $n$ of these coloured cubes to have one Insanity puzzle. The Instant Insanity is just one example constructed by taking $n=4$. There are also Insanity puzzles with $n=5$ that have been in the market as the Hanayama, the Trench Tantalizer, and, the Allies Flag Puzzle (the configurations of these puzzles can be find in \cite{demaine2008lifetime}).  
 
Despite of the huge amount of possibilities to construct an Insanity puzzle, the majority of Insanity puzzles that have been commercialised have the same colouring structure (some of them use four symbols instead of four colours). I share Rik Van Grol's (and O'Bernie's \cite{o1984puzzles}) amazement that people or firms bringing out and commercialising a new Insanity puzzle, do not make an effort to design their own version, but simply replace the original pictures or colours with other figures or colours. This paper is inspired by this amazement.
 
 We have found all possible Insanity puzzles with four, five and six cubes and four five and six colours respectively. We have also classified them by their number of solutions. In doing so we have been able to answer questions like: 

What is the maximum number of solutions that an Insanity puzzle can have? 

Given a number $s$ between this maximum and zero, how many Insanity puzzles have $s$ solutions? 

Once we have fixed one of these Insanity puzzles, is it possible to form a $2 \times 2 \times 1$ prism with monochromatic faces?

We have used combinatorial game theory and graph theory tools to model the Insanity puzzles. The model that we use to analyse Insanity puzzles has been inspired by an algorithm that T.A. Brow presented in \cite{brown1968note}. He has used this algorithm to solve by hand the Instant Insanity puzzle. 

For those interested in creating and designing new Insanity puzzles, we hope that with the results that we present in this paper, it will be possible to construct Insanity puzzles without using the same combinatorial structure of puzzles that are already on the market. Also, these results open the possibility of creating a multilevel Insanity puzzle that can have different puzzles varying with the number of possible solutions.

\section*{Other kind of Insanity Puzzles}

There exists a whole family of variations of Insanity puzzles. Rik Van Grol presents in \cite{demaine2008lifetime} other kind of Insanity puzzles such as On the spot Insanity (that he has designed and presented in the 23rd International Puzzle Collectors Party). 

In the same paper, Rik Van Grol also mentions the Mutando puzzle that consists of four cubes with faces that have four colours but does not fulfill the requirement that all the cubes must have at least one face with each one of the colours. We present in the last section of this paper a detailed description of the Mutando puzzle. We also present a new Insanity puzzle that we called the Mutando of Insanity that is a puzzle analogous to the Mutando but that fulfils the properties of an Insanity puzzle.

In \cite{demaine2013variations}, Demaine et al. defined a whole different family of Insanity puzzles. These puzzles do not necessarily use cubes, they are constructed with regular prisms in such a way that all the coloured faces must apear in the solution of the puzzle.  There exist other great puzzles that consists of coloured cubes and figures that are not considered as Insanity puzzles. We refer the interested reader to a book by P.A. Macmahon \cite{macmahon1921new} were he has written about several coloured figure puzzles. There exists also a generalisation of the Instant Insanity to other Platonic Solids other than the cube \cite{jebasingh2002platonic}. 

These other kind of Insanity puzzles are very interesting and some of them are still in need of a mathematical modelling to analyse them. In this paper we are not going to analyse this other type of Insanity puzzles.

\section*{History of mathematical solutions to the Instant Insanity}

In a personal conversation at G4G 12, David Singmaster gave us references about the first paper that give and analyse a mathematical model of The Tantalizer puzzle (commercialised before the Instant Insanity but with the same colouring structure) \cite{de1947coloured}. In this paper, some (at most four) undergraduate students at Cambridge University, who used to write under the pseudonym F. De Carteblanche, used graph theory to model the Tantalizer and to find the unique solution to the puzzle.

One year after the Instant Insanity was commercialised (in 1967 by Parker Brothers Division), T.A. Brown \cite{brown1968note} used combinatorial number theory to solve the Instant Insanity. 

We have found of mathematical interest that these two mathematical models are equivalent. It is possible to answer questions about possible graph structures with the combinatorial number theory structures that Brown has used and vice versa. In this paper we are not going to talk about this relationship. The mathematical model that we are going to define in the next section has a natural graph theoretical structure and the results that we have found can be settled in the language of graph theory.

As far as we know, the only two mathematical models for the Insanity puzzles that are analysing in this paper are the ones presented in \cite{de1947coloured} and \cite{brown1968note}. Although these are the first papers in which these models were presented, there are other articles and books that talk about Insanity puzzles (\cite{harary1977tantalizer, grecos1971diagrammatic, schwartz1970improved, van1969graph}). Some of these papers give a refinement or a different way of explaining what Brown and De Carteblanche have done. The results that we present in this paper were computed in less than 24 hours for $n=4,5,6$ but it was proved in \cite{robertson1978np} by Eduard Robertson and Ian Munro that the generalisation of Brown's algorithm to all $n \in \N$ is $NP-complete$. 

\section*{Mathematical Modelling of Insanity puzzles}
The definitions that follow were inspired by the algorithm that Brown \cite{brown1968note} used to solve the Instant Insanity. We do not know if the definitions and structures that we are defining in this section have already been studied, but they are a natural way of generalising the matrix structures that Brown has used in \cite{brown1968note}. We believe that by stating the definitions in such a general framework we allow this structures to be useful not only to answer graph theoretical questions (the Insanity puzzles can be transformed in coloured and labeled graphs and be solved with graph theory) but also to be used in more general settings as combinatorial number theory, or algebraic topological structures as simplicial complexes.

From a finite set of prime numbers $\mathcal{P}=\{p_{1},...,p_{n}\}$, we form the set of $k$ products \[\mathcal{P}_{k}=\{p_{i_{1}} \cdot p_{i_{2}}\cdot ... \cdot p_{i_{k}} \mid p_{i_{l}} \in \mathcal{P} \text{ for } l = 1,...,k  \}.\]

For each $m_{1}, m_{2} \in \N$ we define a set containing all matrices with $m_{1}$ rows and $m_{2}$ columns such that each entry of the matrix is an element of $\mathcal{P}_{k}$. We will denote it by $[\mathcal{P}_{k}]^{m_{1}\times m_{2}}$. 

The next definitions will lead us to model the specific structure of the Insanity puzzles with a set of prime numbers  and the set of matrices with prime product entries that we defined above. The set of prime numbers $\mathcal{P}=\{p_{1},...,p_{n}\}$ will represent the set of $n$ colours or labels that we are going to use to colour or mark the faces of the cubes. In the Insanity puzzles, as we have defined them, the number $n$ of labels  is the same as the number of cubes which determine the number $m_{1}=n$ of rows in the matrices $[\mathcal{P}_{k}]^{m_{1}\times m_{2}}$. The labeled structure of each cube will be represented in one row of the matrix in the following way: each entry of a matrix will represent the product of the labels of two opposite parallel faces in a cube and the three opposed faces of each cube will be represented in the three elements in one and only one row. Then, we want to take $k=2$ and $m_{2}=3$. Because our Insanity puzzles have the property that each label must appear in at least one face of each cube, the product of the three elements of each row of a matrix representing an Insanity puzzle must be divisible by $\Pi_{i=1}^{n}p_{n}$. Any matrix of $[\mathcal{P}_{2}]^{n\times 3}$ that has this property will be called a proper matrix. It is easy to show that each possible Insanity puzzle is represented by at least one proper matrix contained in $[\mathcal{P}_{2}]^{n\times 3}$. We are going to denote the subset of all proper matrices by $[\mathcal{P}_{2}]^{n\times 3}_{*}$

It is clear that we can represent an Insanity puzzle with more that one proper matrix but we can define an equivalence relation to be able to have a one to one correspondence between all possible Insanity puzzles and the elements of $[\mathcal{P}_{2}]^{n\times 3}_{*}$. We are going to consider two elements in $[\mathcal{P}_{2}]^{n\times 3}_{*}$ as equivalent if we can get from one to another by permuting the rows or if we can get from one to another by permuting the entries of a row (within the same row).

\begin{eg}[Instant Insanity]
For the Instant Insanity we have $n=4$, we can take $\mathcal{P}=\{2, 3, 5, 7\}$, and then $\mathcal{P}_{2}=\{4, 6, 10, 14, 9, 15, 21, 25, 35, 49\}$. By labelling the faces of the four cubes of the Instant Insanity with the prime numbers in $\mathcal{P}$ instead of colours we can have the following configuration for this puzzle: 

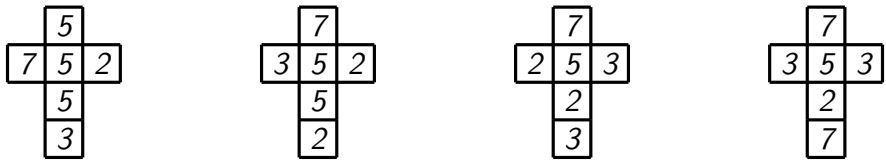
\begin{center}
 \begin{multicols}{4}

\begin{tikzpicture}[thick,scale=.5, ]
   \coordinate (A1) at (0,0);
   \coordinate (A2) at (1, 0);
    \coordinate (A3) at (0, 1);
    \coordinate (A4) at (1, 1);
    \coordinate (B1) at (0, 2);
    \coordinate (B2) at (1, 2);
    \coordinate (B3) at (0, 3);
    \coordinate (B4) at (1, 3);
   \coordinate (C1) at (0, 4);
    \coordinate (C2) at (1, 4);
   \coordinate(D1) at (-1,2);
   \coordinate(D2) at (-1,3);
  \coordinate(D3) at (2,2);
 \coordinate(D4) at (2,3);

   \draw[very thick] (A1) -- (A2);
    \draw[very thick] (A1) -- (A3);
    \draw[very thick] (A2) -- (A4);
    \draw[very thick] (A3) -- (A4);
   \draw[very thick] (A3) -- (B1);
   \draw[very thick] (A4) -- (B2);
   \draw[very thick] (B1) -- (B2);
   \draw[very thick] (B3) -- (C1);
   \draw[very thick] (B4) -- (C2);
   \draw[very thick] (B3) -- (B1);
   \draw[very thick] (B4) -- (B2);
   \draw[very thick] (B3) -- (B4);
 \draw[very thick] (C1) -- (C2);    
  \draw[very thick] (D1) -- (D2); 
 \draw[very thick] (D3) -- (D4);
  \draw[very thick] (B1) -- (D1);
  \draw[very thick] (B2) -- (D3);
  \draw[very thick] (B3) -- (D2);
  \draw[very thick] (B4) -- (D4);

\node at (.5,.5) {3};
\node at (.5,1.5) {5};
\node at (.5,2.5) {5};
\node at (.5,3.5) {5};
\node at (-.5,2.5) {7};
\node at (1.5,2.5) {2};
\end{tikzpicture}

\begin{tikzpicture}[thick,scale=.5,]
   \coordinate (A1) at (0,0);
   \coordinate (A2) at (1, 0);
    \coordinate (A3) at (0, 1);
    \coordinate (A4) at (1, 1);
    \coordinate (B1) at (0, 2);
    \coordinate (B2) at (1, 2);
    \coordinate (B3) at (0, 3);
    \coordinate (B4) at (1, 3);
   \coordinate (C1) at (0, 4);
    \coordinate (C2) at (1, 4);
   \coordinate(D1) at (-1,2);
   \coordinate(D2) at (-1,3);
  \coordinate(D3) at (2,2);
 \coordinate(D4) at (2,3);

   \draw[very thick] (A1) -- (A2);
    \draw[very thick] (A1) -- (A3);
    \draw[very thick] (A2) -- (A4);
    \draw[very thick] (A3) -- (A4);
   \draw[very thick] (A3) -- (B1);
   \draw[very thick] (A4) -- (B2);
   \draw[very thick] (B1) -- (B2);
   \draw[very thick] (B3) -- (C1);
   \draw[very thick] (B4) -- (C2);
   \draw[very thick] (B3) -- (B1);
   \draw[very thick] (B4) -- (B2);
   \draw[very thick] (B3) -- (B4);
 \draw[very thick] (C1) -- (C2);    
  \draw[very thick] (D1) -- (D2); 
 \draw[very thick] (D3) -- (D4);
  \draw[very thick] (B1) -- (D1);
  \draw[very thick] (B2) -- (D3);
  \draw[very thick] (B3) -- (D2);
  \draw[very thick] (B4) -- (D4);

\node at (.5,.5) {2};
\node at (.5,1.5) {5};
\node at (.5,2.5) {5};
\node at (.5,3.5) {7};
\node at (-.5,2.5) {3};
\node at (1.5,2.5) {2};

\end{tikzpicture}

\begin{tikzpicture}[thick,scale=.5,]
   \coordinate (A1) at (0,0);
   \coordinate (A2) at (1, 0);
    \coordinate (A3) at (0, 1);
    \coordinate (A4) at (1, 1);
    \coordinate (B1) at (0, 2);
    \coordinate (B2) at (1, 2);
    \coordinate (B3) at (0, 3);
    \coordinate (B4) at (1, 3);
   \coordinate (C1) at (0, 4);
    \coordinate (C2) at (1, 4);
   \coordinate(D1) at (-1,2);
   \coordinate(D2) at (-1,3);
  \coordinate(D3) at (2,2);
 \coordinate(D4) at (2,3);

   \draw[very thick] (A1) -- (A2) ;
    \draw[very thick] (A1) -- (A3) ;
    \draw[very thick] (A2) -- (A4) ;
    \draw[very thick] (A3) -- (A4);
   \draw[very thick] (A3) -- (B1);
   \draw[very thick] (A4) -- (B2);
   \draw[very thick] (B1) -- (B2);
   \draw[very thick] (B3) -- (C1);
   \draw[very thick] (B4) -- (C2);
   \draw[very thick] (B3) -- (B1);
   \draw[very thick] (B4) -- (B2);
   \draw[very thick] (B3) -- (B4);
 \draw[very thick] (C1) -- (C2);    
  \draw[very thick] (D1) -- (D2); 
 \draw[very thick] (D3) -- (D4);
  \draw[very thick] (B1) -- (D1);
  \draw[very thick] (B2) -- (D3);
  \draw[very thick] (B3) -- (D2);
  \draw[very thick] (B4) -- (D4);

\node at (.5,.5) {3};
\node at (.5,1.5) {2};
\node at (.5,2.5) {5};
\node at (.5,3.5) {7};
\node at (-.5,2.5) {2};
\node at (1.5,2.5) {3};
\end{tikzpicture}

\begin{tikzpicture}[thick,scale=.5, ]
   \coordinate (A1) at (0,0);
   \coordinate (A2) at (1, 0);
    \coordinate (A3) at (0, 1);
    \coordinate (A4) at (1, 1);
    \coordinate (B1) at (0, 2);
    \coordinate (B2) at (1, 2);
    \coordinate (B3) at (0, 3);
    \coordinate (B4) at (1, 3);
   \coordinate (C1) at (0, 4);
    \coordinate (C2) at (1, 4);
   \coordinate(D1) at (-1,2);
   \coordinate(D2) at (-1,3);
  \coordinate(D3) at (2,2);
 \coordinate(D4) at (2,3);

   \draw[very thick] (A1) -- (A2);
    \draw[very thick] (A1) -- (A3);
    \draw[very thick] (A2) -- (A4);
    \draw[very thick] (A3) -- (A4);
   \draw[very thick] (A3) -- (B1);
   \draw[very thick] (A4) -- (B2);
   \draw[very thick] (B1) -- (B2);
   \draw[very thick] (B3) -- (C1);
   \draw[very thick] (B4) -- (C2);
   \draw[very thick] (B3) -- (B1);
   \draw[very thick] (B4) -- (B2);
   \draw[very thick] (B3) -- (B4);
 \draw[very thick] (C1) -- (C2);    
  \draw[very thick] (D1) -- (D2); 
 \draw[very thick] (D3) -- (D4);
  \draw[very thick] (B1) -- (D1);
  \draw[very thick] (B2) -- (D3);
  \draw[very thick] (B3) -- (D2);
  \draw[very thick] (B4) -- (D4);

\node at (.5,.5) {7};
\node at (.5,1.5) {2};
\node at (.5,2.5) {5};
\node at (.5,3.5) {7};
\node at (-.5,2.5) {3};
\node at (1.5,2.5) {3};
\end{tikzpicture}
\end{multicols}
\end{center}
Although we could relabel these four cubes in a different way but still have the configuration of the Instant Insanity by permuting the labels, we are going to consider as different the matrices corresponding to two different ways of labelling the Instant Insanity and we are going to take the one given above as the canonic label of the Instant Insanity.

The next three elements of $[\mathcal{P}_{2}]^{4\times 3}_{*}$ represent the Instant Insanity puzzle and are equivalent because we can get from one to the other by permuting its rows or by permuting the elements of the first and the third rows within the same rows 
\begin{multicols}{3}
\[ \left( \begin{array}{cccc}
14 & 25 & 15 \\
6 & 35 & 10 \\
6 & 14 & 15 \\
9 & 14 & 35 \end{array} \right), \] 

\[ \left( \begin{array}{cccc}
14 & 25 & 15 \\
6 & 14 & 15 \\
6 & 35 & 10 \\
9 & 14 & 35 \end{array} \right),\] 

\[ \left( \begin{array}{cccc}
14 & 15 & 25 \\
6 & 14 & 15 \\
6 & 10 & 35 \\
9 & 14 & 35 \end{array} \right). \] 
\end{multicols}
\end{eg}

The next definitions will allow us to model the solutions of the Insanity puzzles. For an element $A=[a_{(i,j)}] \in [\mathcal{P}_{2}]^{n\times 3}$ and a \textit{magic number} $M$, we define the set of \textit{$M$ partial solutions} for $A$ as $V_{M}^{A}:=\{ \{(1,s_{1}), (2,s_{2}), ..., (n,s_{n})\}\mid  \Pi_{i=1}^{n}a_{(i,s_{i})}=M\}$, where $1 \leq s_{i} \leq 3$ for all $i \in \{1, ..., n\}$. We define two elements of $V_{M}^{A}$ to be \textit{independent} if they have empty intersection. Finally, we define the set of \textit{$l_{M}$-solutions} of $A$ as the set $[S_{M}^{A}]_{l}$ containing the sets of $l$ pairwise independent elements of $V_{M}^{A}$. A solution to an Insanity puzzle consists of four long faces of a tower of $n$ cubes in such a way that each long face of the tower has each one of the $n$ labels. Then, to find all solutions of an Insanity puzzle represented by a matrix $A$, we need two find the set $[S_{M}^{A}]_{2}$ with $M=\Pi_{i=1}^{n}p_{i}^{2}$. We select $l=2$ because there are two opposite long faces in a tower built with $n$ cubes. We choose $M=\Pi_{i=1}^{n}p_{i}^{2}$ because want each label to appear once in each one of the faces of the $n \times 1 \times 1$ tower, which means appearing twice when we consider two opposite faces. 
\begin{eg}[Instant Insanity]
Taking $A$ as the first matrix in Example 1 and $M=2^{2}\cdot 3^{2} \cdot 5^{2} \cdot7^{2}$ it is easy to prove that the only elements of  $ V_{M}^{A}$ are  $\{(1,2), (2,1), (3,2), (4,3)\}$, $\{(1,3), (2,3), (3,2), (4,3)\},$ and $\{(1,3), (2,2), (3,2), (4,2)\} $. Then, the set of $2_{M}$-solutions of $A$ with $M=44,100$ is given by  \[[S_{M}^{A}]_{2}=\{ \{(1,2), (2,1), (3,2), (4,3)\}, \{(1,3), (2,2), (3,2), (4,2)\} \}.\]

We can conclude that the Instant Insanity has only one solution and we can form the solution with the cubes by taking, for example, $\{(1,2), (2,1), (3,2), (4,3)\}$ as the front and back faces of the tower and \{(1,3), (2,2), (3,2), (4,2)\} as the right and left faces of the $4 \times 1 \times $1 tower. 
\end{eg}

In general, for an Insanity puzzle with $n$ labels represented by a proper matrix $A \in [\mathcal{P}_{2}]^{n\times 3}_{*}$ each element of $[S_{M}^{A}]_{2}$ we can form a $n \times 1 \times 1$ tower showing a solution. For this tower if we select its front, back, left, and right faces then the group $S_{n}$ of permutations of four elements and the group $D_{4}$ of the symmetries of the square is acting on the tower and we have $\mid S_{n} \mid \times \mid D_{4} \mid =n! \cdot 8$ different ways of rearranging the cubes in the tower but still have the same solution of the puzzle.

\section*{All possible coloured cubes}
Before finding the number of different Insanity puzzles and the number of solutions that each one has, we need first to know all possible ways for labelling a cube with $n$ labels. With the mathematical model that we have proposed for Insanity puzzles,  what matters is which label is in front of which label in the three different opposite pair of faces of each cube. Because of this, we are going to enumerate bellow all possible sets of three elements that we can form with $\mathcal{P}_{2}$ for $\mathcal{P}=\{2,3,4,5\}$ when $n=4$, $\mathcal{P}=\{2, 3, 4, 5, 7\}$ when $n=5$ and $\mathcal{P}=\{2, 3, 4, 5, 7, 11\}$ when $n=6$. Also we need to guarantee that the product of the three elements is divisible by $\Pi_{i=1}^{n}p_{i}$ because we need each label to be present in each one of the cubes. We are going to denote this set by $[\mathcal{P}_{2}]^{1\times 3}_{n}$. It is an easy task to compute by hand and check with a computer the elements of $[\mathcal{P}_{2}]^{1\times 3}_{n}$ for $n=4,5,6$.
\\

For $\mathcal{P}=\{2, 3, 4, 5\}$ the set $[\mathcal{P}_{2}]^{1\times 3}_{4}$ has $52$ elements. These elements are: 

\begin{center}
\begin{multicols}{4}
\begin{center}
\begin{itemize}
\item[1)] $\{4, 9, 35\}$
\item[2)] $\{4, 25, 21\}$
\item[3)] $\{4, 49, 15\}$
\item[4)] $\{9, 25, 14\}$
\item[5)] $\{9, 49, 10\}$
\item[6)] $\{25, 49, 6\}$
\item[7)] $\{4, 6, 10\}$
\item[8)] $\{4, 10, 21\}$
\item[9)] $\{4, 14, 15\}$
\item[10)] $\{9, 6, 35\}$
\item[11)] $\{9, 15, 14\}$
\item[12)] $\{9, 21, 10\}$
\item[13)] $\{25, 10, 21\}$
\item[14)] $\{25, 15, 14\}$
\item[15)] $\{25, 35, 6\}$
\item[16)] $\{49, 14, 15\}$
\item[17)] $\{49, 21, 10\}$
\item[18)] $\{49, 35, 6\}$
\item[19)] $\{4, 15, 21\}$
\item[20)] $\{4, 15, 35\}$
\item[21)] $\{4, 21, 35\}$
\item[22)] $\{9, 10, 14\}$
\item[23)] $\{9, 10, 35\}$
\item[24)] $\{9, 14, 35\}$
\item[25)] $\{25, 6, 14\}$
\item[26)] $\{25, 6, 21\}$
\item[27)] $\{25, 14, 21\}$
\item[28)] $\{49, 6, 10\}$
\item[29)] $\{49, 6, 15\}$
\item[30)] $\{49, 10, 15\}$
\item[31)] $\{6, 10, 14\}$
\item[32)] $\{6, 15, 21\}$
\item[33)] $\{10, 15, 35\}$
\item[34)] $\{14, 21, 35\}$
\item[35)] $\{6, 6, 35\}$
\item[36)] $\{10, 10, 21\}$
\item[37)] $\{14, 14, 15\}$
\item[38)] $\{15, 15, 14\}$
\item[39)] $\{21, 21, 10\}$
\item[40)] $\{35, 35, 6\}$
\item[29)] $\{49, 6, 15\}$
\item[30)] $\{49, 10, 15\}$
\item[31)] $\{6, 10, 14\}$
\item[32)] $\{6, 15, 21\}$
\item[33)] $\{10, 15, 35\}$
\item[34)] $\{14, 21, 35\}$
\item[35)] $\{6, 6, 35\}$
\item[36)] $\{10, 10, 21\}$
\item[37)] $\{14, 14, 15\}$
\item[38)] $\{15, 15, 14\}$
\item[39)] $\{21, 21, 10\}$
\item[40)] $\{35, 35, 6\}$
\item[41)] $\{6, 10, 21\}$
\item[42)] $\{6, 10, 35\}$
\item[43)] $\{6, 14, 15\}$
\item[44)] $\{6, 14, 35\}$
\item[45)] $\{6, 15, 35\}$
\item[46)] $\{6, 21, 35\}$
\item[47)] $\{10, 14, 15\}$
\item[48)] $\{10, 14, 21\}$
\item[49)] $\{10, 15, 21\}$
\item[50)] $\{10, 35, 21\}$
\item[51)] $\{15, 35, 14\}$
\item[52)] $\{14, 21, 15\}$
\end{itemize}
\end{center}
\end{multicols}
\end{center}

For $\mathcal{P}=\{2, 3, 4, 5, 7\}$ the set $[\mathcal{P}_{2}]^{1\times 3}_{5}$ has $45$ elements. These elements are:
 
\begin{center}
\begin{multicols}{4}
\begin{center}
\begin{itemize}
\item[1)] $\{4, 15, 77\}$
\item[2)] $\{4, 21, 55\}$
\item[3)] $\{4, 33, 35\}$
\item[4)] $\{9, 10, 77\}$
\item[5)] $\{9, 14, 55\}$
\item[6)] $\{9, 22, 35\}$
\item[7)] $\{25, 6, 77\}$
\item[8)] $\{25, 14, 15\}$
\item[9)] $\{25, 22, 21\}$
\item[10)] $\{49, 6, 55\}$
\item[11)] $\{49, 10, 33\}$
\item[12)] $\{49, 22, 15\}$
\item[13)] $\{121, 6, 35\}$
\item[14)] $\{121, 10, 21\}$
\item[15)] $\{121, 14, 15\}$
\item[16)] $\{6, 10, 77\}$
\item[17)] $\{6, 14, 55\}$
\item[18)] $\{6, 22, 35\}$
\item[19)] $\{10, 14, 33\}$
\item[20)] $\{10, 22, 21\}$
\item[21)] $\{14, 22, 15\}$
\item[22)] $\{6, 15, 77\}$
\item[23)] $\{6, 21, 55\}$
\item[24)] $\{6, 33, 35\}$
\item[25)] $\{15, 21, 22\}$
\item[26)] $\{15, 33, 14\}$
\item[27)] $\{21, 33, 10\}$
\item[28)] $\{10, 15, 77\}$
\item[29)] $\{10, 35, 33\}$
\item[30)] $\{10, 55, 21\}$
\item[31)] $\{15, 35, 22\}$
\item[32)] $\{15, 55, 14\}$
\item[33)] $\{35, 55, 10\}$
\item[34)] $\{14, 21, 35\}$
\item[35)] $\{14, 35, 33\}$
\item[36)] $\{14, 77, 15\}$
\item[37)] $\{21, 35, 22\}$
\item[38)] $\{21, 77, 10\}$
\item[39)] $\{35, 77, 6\}$
\item[40)] $\{22, 33, 35\}$
\item[41)] $\{22, 55, 21\}$
\item[42)] $\{22, 77, 15\}$
\item[43)] $\{33, 55, 14\}$
\item[44)] $\{33, 77, 10\}$
\item[45)] $\{55, 77, 6\}$
\end{itemize}
\end{center}
\end{multicols}
\end{center}

For $\mathcal{P}=\{2, 3, 4, 5, 11, 13\}$ the set $[\mathcal{P}_{2}]^{1\times 3}_{6}$ has $15$ elements. These elements are:
 
\begin{center}
\begin{multicols}{4}
\begin{center}
\begin{itemize}
\item[1)] $\{6, 35, 143\}$
\item[2)] $\{6, 55, 91\}$
\item[3)] $\{6, 65, 77\}$
\item[4)] $\{10, 21, 143\}$
\item[5)] $\{10, 33, 91\}$
\item[6)] $\{10, 39, 77\}$
\item[7)] $\{14, 15, 143\}$
\item[8)] $\{14, 33, 65\}$
\item[9)] $\{14, 39, 55\}$
\item[10)] $\{22, 15, 91\}$
\item[11)] $\{22, 21, 65\}$
\item[12)] $\{22, 39, 35\}$
\item[13)] $\{26, 15, 77\}$
\item[14)] $\{26, 35, 55\}$
\item[15)] $\{26, 33, 35\}$
\end{itemize}
\end{center}
\end{multicols}
\end{center}

\section*{Solutions to all Insanity puzzles}
The purpose of this section is to analyse all possible Insanity puzzles for $n \in \{4,5,6\}$. We were able to find the number of solutions to all Insanity puzzles. In particular, for $n=4$ we can know how many Insanity puzzles can be designed that have a different abstract structure with only one solution. This gives us a whole set of different Insanity puzzles with unique solution that is different from the Instant Insanity puzzle. To form a matrix $A$ representing an Insanity puzzle, we take four distinct elements of $[\mathcal{P}_{2}]^{1\times 3}_{n}$ as rows to form an $n \times 3$ matrix. Then we calculate the set of solutions  $[S_{M}^{A}]_{2}$ and its cardinality will tell us how many different solutions has the Insanity puzzle represented by $A$. Based on the mathematical model that we have constructed for the Insanity puzzles, we have implemented an algorithm to calculate all Insanity puzzles. We have found the next results.


For $n=4$ there exists an Insanity puzzle with $72$ solutions and there are no Insanity puzzles with more than $72$ solutions.
It is not true that for all $0 \leq m \leq 72$ there exists an Insanity puzzle with $m$ solutions. For $m= 13, 15, 17, 19, 22, 23, 25, 26, 27, 29, 30, 31, 32, 33, 34, 35, 37, 38, 39, 40$ $41, 42, 43, 44, 45, 46, 47, 49, 50, 51, 52, 53, 54, 55, 56, 57, 58, 59, 60, 61, 62, 63, 64, 65, 66, 67, 68, 69,$ $70$, or $m=71$, there is no Insanity puzzle with $m$ distinct solutions.
 
The next Insanity puzzle with four cubes has only one solution and its configuration is different from the Instant Insanity.

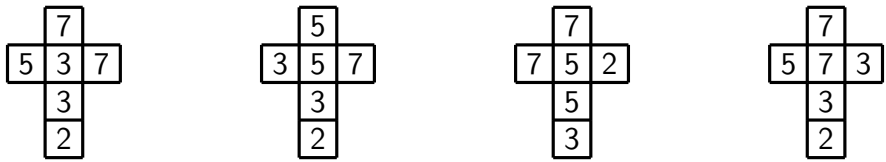
\begin{center}
 \begin{multicols}{4}

\begin{tikzpicture}[thick,scale=.5, ]
   \coordinate (A1) at (0,0);
   \coordinate (A2) at (1, 0);
    \coordinate (A3) at (0, 1);
    \coordinate (A4) at (1, 1);
    \coordinate (B1) at (0, 2);
    \coordinate (B2) at (1, 2);
    \coordinate (B3) at (0, 3);
    \coordinate (B4) at (1, 3);
   \coordinate (C1) at (0, 4);
    \coordinate (C2) at (1, 4);
   \coordinate(D1) at (-1,2);
   \coordinate(D2) at (-1,3);
  \coordinate(D3) at (2,2);
 \coordinate(D4) at (2,3);

   \draw[very thick] (A1) -- (A2);
    \draw[very thick] (A1) -- (A3);
    \draw[very thick] (A2) -- (A4);
    \draw[very thick] (A3) -- (A4);
   \draw[very thick] (A3) -- (B1);
   \draw[very thick] (A4) -- (B2);
   \draw[very thick] (B1) -- (B2);
   \draw[very thick] (B3) -- (C1);
   \draw[very thick] (B4) -- (C2);
   \draw[very thick] (B3) -- (B1);
   \draw[very thick] (B4) -- (B2);
   \draw[very thick] (B3) -- (B4);
 \draw[very thick] (C1) -- (C2);    
  \draw[very thick] (D1) -- (D2); 
 \draw[very thick] (D3) -- (D4);
  \draw[very thick] (B1) -- (D1);
  \draw[very thick] (B2) -- (D3);
  \draw[very thick] (B3) -- (D2);
  \draw[very thick] (B4) -- (D4);

\node at (.5,.5) {2};
\node at (.5,1.5) {3};
\node at (.5,2.5) {3};
\node at (.5,3.5) {7};
\node at (-.5,2.5) {5};
\node at (1.5,2.5) {7};
\end{tikzpicture}

\begin{tikzpicture}[thick,scale=.5,]
   \coordinate (A1) at (0,0);
   \coordinate (A2) at (1, 0);
    \coordinate (A3) at (0, 1);
    \coordinate (A4) at (1, 1);
    \coordinate (B1) at (0, 2);
    \coordinate (B2) at (1, 2);
    \coordinate (B3) at (0, 3);
    \coordinate (B4) at (1, 3);
   \coordinate (C1) at (0, 4);
    \coordinate (C2) at (1, 4);
   \coordinate(D1) at (-1,2);
   \coordinate(D2) at (-1,3);
  \coordinate(D3) at (2,2);
 \coordinate(D4) at (2,3);

   \draw[very thick] (A1) -- (A2);
    \draw[very thick] (A1) -- (A3);
    \draw[very thick] (A2) -- (A4);
    \draw[very thick] (A3) -- (A4);
   \draw[very thick] (A3) -- (B1);
   \draw[very thick] (A4) -- (B2);
   \draw[very thick] (B1) -- (B2);
   \draw[very thick] (B3) -- (C1);
   \draw[very thick] (B4) -- (C2);
   \draw[very thick] (B3) -- (B1);
   \draw[very thick] (B4) -- (B2);
   \draw[very thick] (B3) -- (B4);
 \draw[very thick] (C1) -- (C2);    
  \draw[very thick] (D1) -- (D2); 
 \draw[very thick] (D3) -- (D4);
  \draw[very thick] (B1) -- (D1);
  \draw[very thick] (B2) -- (D3);
  \draw[very thick] (B3) -- (D2);
  \draw[very thick] (B4) -- (D4);

\node at (.5,.5) {2};
\node at (.5,1.5) {3};
\node at (.5,2.5) {5};
\node at (.5,3.5) {5};
\node at (-.5,2.5) {3};
\node at (1.5,2.5) {7};

\end{tikzpicture}

\begin{tikzpicture}[thick,scale=.5,]
   \coordinate (A1) at (0,0);
   \coordinate (A2) at (1, 0);
    \coordinate (A3) at (0, 1);
    \coordinate (A4) at (1, 1);
    \coordinate (B1) at (0, 2);
    \coordinate (B2) at (1, 2);
    \coordinate (B3) at (0, 3);
    \coordinate (B4) at (1, 3);
   \coordinate (C1) at (0, 4);
    \coordinate (C2) at (1, 4);
   \coordinate(D1) at (-1,2);
   \coordinate(D2) at (-1,3);
  \coordinate(D3) at (2,2);
 \coordinate(D4) at (2,3);

   \draw[very thick] (A1) -- (A2) ;
    \draw[very thick] (A1) -- (A3) ;
    \draw[very thick] (A2) -- (A4) ;
    \draw[very thick] (A3) -- (A4);
   \draw[very thick] (A3) -- (B1);
   \draw[very thick] (A4) -- (B2);
   \draw[very thick] (B1) -- (B2);
   \draw[very thick] (B3) -- (C1);
   \draw[very thick] (B4) -- (C2);
   \draw[very thick] (B3) -- (B1);
   \draw[very thick] (B4) -- (B2);
   \draw[very thick] (B3) -- (B4);
 \draw[very thick] (C1) -- (C2);    
  \draw[very thick] (D1) -- (D2); 
 \draw[very thick] (D3) -- (D4);
  \draw[very thick] (B1) -- (D1);
  \draw[very thick] (B2) -- (D3);
  \draw[very thick] (B3) -- (D2);
  \draw[very thick] (B4) -- (D4);

\node at (.5,.5) {3};
\node at (.5,1.5) {5};
\node at (.5,2.5) {5};
\node at (.5,3.5) {7};
\node at (-.5,2.5) {7};
\node at (1.5,2.5) {2};
\end{tikzpicture}

\begin{tikzpicture}[thick,scale=.5, ]
   \coordinate (A1) at (0,0);
   \coordinate (A2) at (1, 0);
    \coordinate (A3) at (0, 1);
    \coordinate (A4) at (1, 1);
    \coordinate (B1) at (0, 2);
    \coordinate (B2) at (1, 2);
    \coordinate (B3) at (0, 3);
    \coordinate (B4) at (1, 3);
   \coordinate (C1) at (0, 4);
    \coordinate (C2) at (1, 4);
   \coordinate(D1) at (-1,2);
   \coordinate(D2) at (-1,3);
  \coordinate(D3) at (2,2);
 \coordinate(D4) at (2,3);

   \draw[very thick] (A1) -- (A2);
    \draw[very thick] (A1) -- (A3);
    \draw[very thick] (A2) -- (A4);
    \draw[very thick] (A3) -- (A4);
   \draw[very thick] (A3) -- (B1);
   \draw[very thick] (A4) -- (B2);
   \draw[very thick] (B1) -- (B2);
   \draw[very thick] (B3) -- (C1);
   \draw[very thick] (B4) -- (C2);
   \draw[very thick] (B3) -- (B1);
   \draw[very thick] (B4) -- (B2);
   \draw[very thick] (B3) -- (B4);
 \draw[very thick] (C1) -- (C2);    
  \draw[very thick] (D1) -- (D2); 
 \draw[very thick] (D3) -- (D4);
  \draw[very thick] (B1) -- (D1);
  \draw[very thick] (B2) -- (D3);
  \draw[very thick] (B3) -- (D2);
  \draw[very thick] (B4) -- (D4);

\node at (.5,.5) {2};
\node at (.5,1.5) {3};
\node at (.5,2.5) {7};
\node at (.5,3.5) {7};
\node at (-.5,2.5) {5};
\node at (1.5,2.5) {3};
\end{tikzpicture}
\end{multicols}
\end{center}

The next Insanity puzzle with four cubes has $72$ solutions.
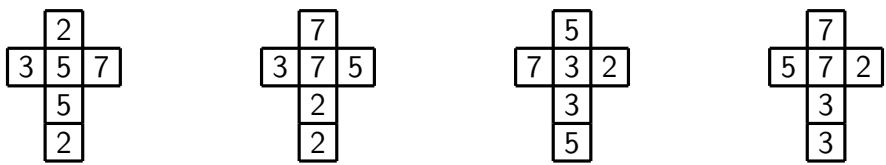
\begin{center}
 \begin{multicols}{4}
\begin{tikzpicture}[thick,scale=.5, ]
   \coordinate (A1) at (0,0);
   \coordinate (A2) at (1, 0);
    \coordinate (A3) at (0, 1);
    \coordinate (A4) at (1, 1);
    \coordinate (B1) at (0, 2);
    \coordinate (B2) at (1, 2);
    \coordinate (B3) at (0, 3);
    \coordinate (B4) at (1, 3);
   \coordinate (C1) at (0, 4);
    \coordinate (C2) at (1, 4);
   \coordinate(D1) at (-1,2);
   \coordinate(D2) at (-1,3);
  \coordinate(D3) at (2,2);
 \coordinate(D4) at (2,3);

   \draw[very thick] (A1) -- (A2);
    \draw[very thick] (A1) -- (A3);
    \draw[very thick] (A2) -- (A4);
    \draw[very thick] (A3) -- (A4);
   \draw[very thick] (A3) -- (B1);
   \draw[very thick] (A4) -- (B2);
   \draw[very thick] (B1) -- (B2);
   \draw[very thick] (B3) -- (C1);
   \draw[very thick] (B4) -- (C2);
   \draw[very thick] (B3) -- (B1);
   \draw[very thick] (B4) -- (B2);
   \draw[very thick] (B3) -- (B4);
 \draw[very thick] (C1) -- (C2);    
  \draw[very thick] (D1) -- (D2); 
 \draw[very thick] (D3) -- (D4);
  \draw[very thick] (B1) -- (D1);
  \draw[very thick] (B2) -- (D3);
  \draw[very thick] (B3) -- (D2);
  \draw[very thick] (B4) -- (D4);

\node at (.5,.5) {2};
\node at (.5,1.5) {5};
\node at (.5,2.5) {5};
\node at (.5,3.5) {2};
\node at (-.5,2.5) {3};
\node at (1.5,2.5) {7};
\end{tikzpicture}

\begin{tikzpicture}[thick,scale=.5,]
   \coordinate (A1) at (0,0);
   \coordinate (A2) at (1, 0);
    \coordinate (A3) at (0, 1);
    \coordinate (A4) at (1, 1);
    \coordinate (B1) at (0, 2);
    \coordinate (B2) at (1, 2);
    \coordinate (B3) at (0, 3);
    \coordinate (B4) at (1, 3);
   \coordinate (C1) at (0, 4);
    \coordinate (C2) at (1, 4);
   \coordinate(D1) at (-1,2);
   \coordinate(D2) at (-1,3);
  \coordinate(D3) at (2,2);
 \coordinate(D4) at (2,3);

   \draw[very thick] (A1) -- (A2);
    \draw[very thick] (A1) -- (A3);
    \draw[very thick] (A2) -- (A4);
    \draw[very thick] (A3) -- (A4);
   \draw[very thick] (A3) -- (B1);
   \draw[very thick] (A4) -- (B2);
   \draw[very thick] (B1) -- (B2);
   \draw[very thick] (B3) -- (C1);
   \draw[very thick] (B4) -- (C2);
   \draw[very thick] (B3) -- (B1);
   \draw[very thick] (B4) -- (B2);
   \draw[very thick] (B3) -- (B4);
 \draw[very thick] (C1) -- (C2);    
  \draw[very thick] (D1) -- (D2); 
 \draw[very thick] (D3) -- (D4);
  \draw[very thick] (B1) -- (D1);
  \draw[very thick] (B2) -- (D3);
  \draw[very thick] (B3) -- (D2);
  \draw[very thick] (B4) -- (D4);

\node at (.5,.5) {2};
\node at (.5,1.5) {2};
\node at (.5,2.5) {7};
\node at (.5,3.5) {7};
\node at (-.5,2.5) {3};
\node at (1.5,2.5) {5};

\end{tikzpicture}

\begin{tikzpicture}[thick,scale=.5,]
   \coordinate (A1) at (0,0);
   \coordinate (A2) at (1, 0);
    \coordinate (A3) at (0, 1);
    \coordinate (A4) at (1, 1);
    \coordinate (B1) at (0, 2);
    \coordinate (B2) at (1, 2);
    \coordinate (B3) at (0, 3);
    \coordinate (B4) at (1, 3);
   \coordinate (C1) at (0, 4);
    \coordinate (C2) at (1, 4);
   \coordinate(D1) at (-1,2);
   \coordinate(D2) at (-1,3);
  \coordinate(D3) at (2,2);
 \coordinate(D4) at (2,3);

   \draw[very thick] (A1) -- (A2) ;
    \draw[very thick] (A1) -- (A3) ;
    \draw[very thick] (A2) -- (A4) ;
    \draw[very thick] (A3) -- (A4);
   \draw[very thick] (A3) -- (B1);
   \draw[very thick] (A4) -- (B2);
   \draw[very thick] (B1) -- (B2);
   \draw[very thick] (B3) -- (C1);
   \draw[very thick] (B4) -- (C2);
   \draw[very thick] (B3) -- (B1);
   \draw[very thick] (B4) -- (B2);
   \draw[very thick] (B3) -- (B4);
 \draw[very thick] (C1) -- (C2);    
  \draw[very thick] (D1) -- (D2); 
 \draw[very thick] (D3) -- (D4);
  \draw[very thick] (B1) -- (D1);
  \draw[very thick] (B2) -- (D3);
  \draw[very thick] (B3) -- (D2);
  \draw[very thick] (B4) -- (D4);

\node at (.5,.5) {5};
\node at (.5,1.5) {3};
\node at (.5,2.5) {3};
\node at (.5,3.5) {5};
\node at (-.5,2.5) {7};
\node at (1.5,2.5) {2};
\end{tikzpicture}

\begin{tikzpicture}[thick,scale=.5, ]
   \coordinate (A1) at (0,0);
   \coordinate (A2) at (1, 0);
    \coordinate (A3) at (0, 1);
    \coordinate (A4) at (1, 1);
    \coordinate (B1) at (0, 2);
    \coordinate (B2) at (1, 2);
    \coordinate (B3) at (0, 3);
    \coordinate (B4) at (1, 3);
   \coordinate (C1) at (0, 4);
    \coordinate (C2) at (1, 4);
   \coordinate(D1) at (-1,2);
   \coordinate(D2) at (-1,3);
  \coordinate(D3) at (2,2);
 \coordinate(D4) at (2,3);

   \draw[very thick] (A1) -- (A2);
    \draw[very thick] (A1) -- (A3);
    \draw[very thick] (A2) -- (A4);
    \draw[very thick] (A3) -- (A4);
   \draw[very thick] (A3) -- (B1);
   \draw[very thick] (A4) -- (B2);
   \draw[very thick] (B1) -- (B2);
   \draw[very thick] (B3) -- (C1);
   \draw[very thick] (B4) -- (C2);
   \draw[very thick] (B3) -- (B1);
   \draw[very thick] (B4) -- (B2);
   \draw[very thick] (B3) -- (B4);
 \draw[very thick] (C1) -- (C2);    
  \draw[very thick] (D1) -- (D2); 
 \draw[very thick] (D3) -- (D4);
  \draw[very thick] (B1) -- (D1);
  \draw[very thick] (B2) -- (D3);
  \draw[very thick] (B3) -- (D2);
  \draw[very thick] (B4) -- (D4);

\node at (.5,.5) {3};
\node at (.5,1.5) {3};
\node at (.5,2.5) {7};
\node at (.5,3.5) {7};
\node at (-.5,2.5) {5};
\node at (1.5,2.5) {2};
\end{tikzpicture}
\end{multicols}
\end{center}

For $n=5$ there exists an Insanity puzzle with $18$ solutions and there are no Insanity puzzles with more than $18$ solutions. It is not true that for all $0 \leq m \leq 18$ there exists a Insanity puzzle with $m$ solutions. For $m=14$ and $m=15$ there is no Insanity puzzle with $5$ cubes and $m$ distinct solutions.

The next Insanity puzzle with five cubes has only one solution.
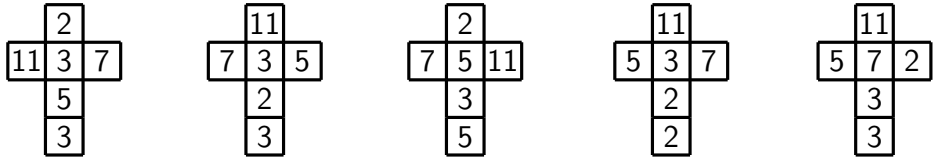
\begin{center}
 \begin{multicols}{5}

\begin{tikzpicture}[thick,scale=.5, ]
   \coordinate (A1) at (0,0);
   \coordinate (A2) at (1, 0);
    \coordinate (A3) at (0, 1);
    \coordinate (A4) at (1, 1);
    \coordinate (B1) at (0, 2);
    \coordinate (B2) at (1, 2);
    \coordinate (B3) at (0, 3);
    \coordinate (B4) at (1, 3);
   \coordinate (C1) at (0, 4);
    \coordinate (C2) at (1, 4);
   \coordinate(D1) at (-1,2);
   \coordinate(D2) at (-1,3);
  \coordinate(D3) at (2,2);
 \coordinate(D4) at (2,3);

   \draw[very thick] (A1) -- (A2);
    \draw[very thick] (A1) -- (A3);
    \draw[very thick] (A2) -- (A4);
    \draw[very thick] (A3) -- (A4);
   \draw[very thick] (A3) -- (B1);
   \draw[very thick] (A4) -- (B2);
   \draw[very thick] (B1) -- (B2);
   \draw[very thick] (B3) -- (C1);
   \draw[very thick] (B4) -- (C2);
   \draw[very thick] (B3) -- (B1);
   \draw[very thick] (B4) -- (B2);
   \draw[very thick] (B3) -- (B4);
 \draw[very thick] (C1) -- (C2);    
  \draw[very thick] (D1) -- (D2); 
 \draw[very thick] (D3) -- (D4);
  \draw[very thick] (B1) -- (D1);
  \draw[very thick] (B2) -- (D3);
  \draw[very thick] (B3) -- (D2);
  \draw[very thick] (B4) -- (D4);

\node at (.5,.5) {3};
\node at (.5,1.5) {5};
\node at (.5,2.5) {3};
\node at (.5,3.5) {2};
\node at (-.5,2.5) {11};
\node at (1.5,2.5) {7};
\end{tikzpicture}

\begin{tikzpicture}[thick,scale=.5,]
   \coordinate (A1) at (0,0);
   \coordinate (A2) at (1, 0);
    \coordinate (A3) at (0, 1);
    \coordinate (A4) at (1, 1);
    \coordinate (B1) at (0, 2);
    \coordinate (B2) at (1, 2);
    \coordinate (B3) at (0, 3);
    \coordinate (B4) at (1, 3);
   \coordinate (C1) at (0, 4);
    \coordinate (C2) at (1, 4);
   \coordinate(D1) at (-1,2);
   \coordinate(D2) at (-1,3);
  \coordinate(D3) at (2,2);
 \coordinate(D4) at (2,3);

   \draw[very thick] (A1) -- (A2);
    \draw[very thick] (A1) -- (A3);
    \draw[very thick] (A2) -- (A4);
    \draw[very thick] (A3) -- (A4);
   \draw[very thick] (A3) -- (B1);
   \draw[very thick] (A4) -- (B2);
   \draw[very thick] (B1) -- (B2);
   \draw[very thick] (B3) -- (C1);
   \draw[very thick] (B4) -- (C2);
   \draw[very thick] (B3) -- (B1);
   \draw[very thick] (B4) -- (B2);
   \draw[very thick] (B3) -- (B4);
 \draw[very thick] (C1) -- (C2);    
  \draw[very thick] (D1) -- (D2); 
 \draw[very thick] (D3) -- (D4);
  \draw[very thick] (B1) -- (D1);
  \draw[very thick] (B2) -- (D3);
  \draw[very thick] (B3) -- (D2);
  \draw[very thick] (B4) -- (D4);

\node at (.5,.5) {3};
\node at (.5,1.5) {2};
\node at (.5,2.5) {3};
\node at (.5,3.5) {11};
\node at (-.5,2.5) {7};
\node at (1.5,2.5) {5};

\end{tikzpicture}

\begin{tikzpicture}[thick,scale=.5,]
   \coordinate (A1) at (0,0);
   \coordinate (A2) at (1, 0);
    \coordinate (A3) at (0, 1);
    \coordinate (A4) at (1, 1);
    \coordinate (B1) at (0, 2);
    \coordinate (B2) at (1, 2);
    \coordinate (B3) at (0, 3);
    \coordinate (B4) at (1, 3);
   \coordinate (C1) at (0, 4);
    \coordinate (C2) at (1, 4);
   \coordinate(D1) at (-1,2);
   \coordinate(D2) at (-1,3);
  \coordinate(D3) at (2,2);
 \coordinate(D4) at (2,3);

   \draw[very thick] (A1) -- (A2) ;
    \draw[very thick] (A1) -- (A3) ;
    \draw[very thick] (A2) -- (A4) ;
    \draw[very thick] (A3) -- (A4);
   \draw[very thick] (A3) -- (B1);
   \draw[very thick] (A4) -- (B2);
   \draw[very thick] (B1) -- (B2);
   \draw[very thick] (B3) -- (C1);
   \draw[very thick] (B4) -- (C2);
   \draw[very thick] (B3) -- (B1);
   \draw[very thick] (B4) -- (B2);
   \draw[very thick] (B3) -- (B4);
 \draw[very thick] (C1) -- (C2);    
  \draw[very thick] (D1) -- (D2); 
 \draw[very thick] (D3) -- (D4);
  \draw[very thick] (B1) -- (D1);
  \draw[very thick] (B2) -- (D3);
  \draw[very thick] (B3) -- (D2);
  \draw[very thick] (B4) -- (D4);

\node at (.5,.5) {5};
\node at (.5,1.5) {3};
\node at (.5,2.5) {5};
\node at (.5,3.5) {2};
\node at (-.5,2.5) {7};
\node at (1.5,2.5) {11};
\end{tikzpicture}

\begin{tikzpicture}[thick,scale=.5, ]
   \coordinate (A1) at (0,0);
   \coordinate (A2) at (1, 0);
    \coordinate (A3) at (0, 1);
    \coordinate (A4) at (1, 1);
    \coordinate (B1) at (0, 2);
    \coordinate (B2) at (1, 2);
    \coordinate (B3) at (0, 3);
    \coordinate (B4) at (1, 3);
   \coordinate (C1) at (0, 4);
    \coordinate (C2) at (1, 4);
   \coordinate(D1) at (-1,2);
   \coordinate(D2) at (-1,3);
  \coordinate(D3) at (2,2);
 \coordinate(D4) at (2,3);

   \draw[very thick] (A1) -- (A2);
    \draw[very thick] (A1) -- (A3);
    \draw[very thick] (A2) -- (A4);
    \draw[very thick] (A3) -- (A4);
   \draw[very thick] (A3) -- (B1);
   \draw[very thick] (A4) -- (B2);
   \draw[very thick] (B1) -- (B2);
   \draw[very thick] (B3) -- (C1);
   \draw[very thick] (B4) -- (C2);
   \draw[very thick] (B3) -- (B1);
   \draw[very thick] (B4) -- (B2);
   \draw[very thick] (B3) -- (B4);
 \draw[very thick] (C1) -- (C2);    
  \draw[very thick] (D1) -- (D2); 
 \draw[very thick] (D3) -- (D4);
  \draw[very thick] (B1) -- (D1);
  \draw[very thick] (B2) -- (D3);
  \draw[very thick] (B3) -- (D2);
  \draw[very thick] (B4) -- (D4);

\node at (.5,.5) {2};
\node at (.5,1.5) {2};
\node at (.5,2.5) {3};
\node at (.5,3.5) {11};
\node at (-.5,2.5) {5};
\node at (1.5,2.5) {7};
\end{tikzpicture}

\begin{tikzpicture}[thick,scale=.5, ]
   \coordinate (A1) at (0,0);
   \coordinate (A2) at (1, 0);
    \coordinate (A3) at (0, 1);
    \coordinate (A4) at (1, 1);
    \coordinate (B1) at (0, 2);
    \coordinate (B2) at (1, 2);
    \coordinate (B3) at (0, 3);
    \coordinate (B4) at (1, 3);
   \coordinate (C1) at (0, 4);
    \coordinate (C2) at (1, 4);
   \coordinate(D1) at (-1,2);
   \coordinate(D2) at (-1,3);
  \coordinate(D3) at (2,2);
 \coordinate(D4) at (2,3);

   \draw[very thick] (A1) -- (A2);
    \draw[very thick] (A1) -- (A3);
    \draw[very thick] (A2) -- (A4);
    \draw[very thick] (A3) -- (A4);
   \draw[very thick] (A3) -- (B1);
   \draw[very thick] (A4) -- (B2);
   \draw[very thick] (B1) -- (B2);
   \draw[very thick] (B3) -- (C1);
   \draw[very thick] (B4) -- (C2);
   \draw[very thick] (B3) -- (B1);
   \draw[very thick] (B4) -- (B2);
   \draw[very thick] (B3) -- (B4);
 \draw[very thick] (C1) -- (C2);    
  \draw[very thick] (D1) -- (D2); 
 \draw[very thick] (D3) -- (D4);
  \draw[very thick] (B1) -- (D1);
  \draw[very thick] (B2) -- (D3);
  \draw[very thick] (B3) -- (D2);
  \draw[very thick] (B4) -- (D4);

\node at (.5,.5) {3};
\node at (.5,1.5) {3};
\node at (.5,2.5) {7};
\node at (.5,3.5) {11};
\node at (-.5,2.5) {5};
\node at (1.5,2.5) {2};
\end{tikzpicture}

\end{multicols}
\end{center}
The next Insanity puzzle with five cubes has $18$ solutions.
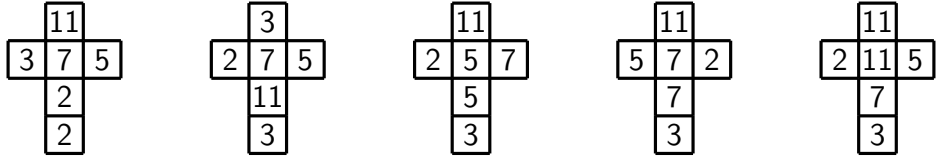
\begin{center}
 \begin{multicols}{5}

\begin{tikzpicture}[thick,scale=.5, ]
   \coordinate (A1) at (0,0);
   \coordinate (A2) at (1, 0);
    \coordinate (A3) at (0, 1);
    \coordinate (A4) at (1, 1);
    \coordinate (B1) at (0, 2);
    \coordinate (B2) at (1, 2);
    \coordinate (B3) at (0, 3);
    \coordinate (B4) at (1, 3);
   \coordinate (C1) at (0, 4);
    \coordinate (C2) at (1, 4);
   \coordinate(D1) at (-1,2);
   \coordinate(D2) at (-1,3);
  \coordinate(D3) at (2,2);
 \coordinate(D4) at (2,3);

   \draw[very thick] (A1) -- (A2);
    \draw[very thick] (A1) -- (A3);
    \draw[very thick] (A2) -- (A4);
    \draw[very thick] (A3) -- (A4);
   \draw[very thick] (A3) -- (B1);
   \draw[very thick] (A4) -- (B2);
   \draw[very thick] (B1) -- (B2);
   \draw[very thick] (B3) -- (C1);
   \draw[very thick] (B4) -- (C2);
   \draw[very thick] (B3) -- (B1);
   \draw[very thick] (B4) -- (B2);
   \draw[very thick] (B3) -- (B4);
 \draw[very thick] (C1) -- (C2);    
  \draw[very thick] (D1) -- (D2); 
 \draw[very thick] (D3) -- (D4);
  \draw[very thick] (B1) -- (D1);
  \draw[very thick] (B2) -- (D3);
  \draw[very thick] (B3) -- (D2);
  \draw[very thick] (B4) -- (D4);

\node at (.5,.5) {2};
\node at (.5,1.5) {2};
\node at (.5,2.5) {7};
\node at (.5,3.5) {11};
\node at (-.5,2.5) {3};
\node at (1.5,2.5) {5};
\end{tikzpicture}

\begin{tikzpicture}[thick,scale=.5,]
   \coordinate (A1) at (0,0);
   \coordinate (A2) at (1, 0);
    \coordinate (A3) at (0, 1);
    \coordinate (A4) at (1, 1);
    \coordinate (B1) at (0, 2);
    \coordinate (B2) at (1, 2);
    \coordinate (B3) at (0, 3);
    \coordinate (B4) at (1, 3);
   \coordinate (C1) at (0, 4);
    \coordinate (C2) at (1, 4);
   \coordinate(D1) at (-1,2);
   \coordinate(D2) at (-1,3);
  \coordinate(D3) at (2,2);
 \coordinate(D4) at (2,3);

   \draw[very thick] (A1) -- (A2);
    \draw[very thick] (A1) -- (A3);
    \draw[very thick] (A2) -- (A4);
    \draw[very thick] (A3) -- (A4);
   \draw[very thick] (A3) -- (B1);
   \draw[very thick] (A4) -- (B2);
   \draw[very thick] (B1) -- (B2);
   \draw[very thick] (B3) -- (C1);
   \draw[very thick] (B4) -- (C2);
   \draw[very thick] (B3) -- (B1);
   \draw[very thick] (B4) -- (B2);
   \draw[very thick] (B3) -- (B4);
 \draw[very thick] (C1) -- (C2);    
  \draw[very thick] (D1) -- (D2); 
 \draw[very thick] (D3) -- (D4);
  \draw[very thick] (B1) -- (D1);
  \draw[very thick] (B2) -- (D3);
  \draw[very thick] (B3) -- (D2);
  \draw[very thick] (B4) -- (D4);

\node at (.5,.5) {3};
\node at (.5,1.5) {11};
\node at (.5,2.5) {7};
\node at (.5,3.5) {3};
\node at (-.5,2.5) {2};
\node at (1.5,2.5) {5};

\end{tikzpicture}

\begin{tikzpicture}[thick,scale=.5,]
   \coordinate (A1) at (0,0);
   \coordinate (A2) at (1, 0);
    \coordinate (A3) at (0, 1);
    \coordinate (A4) at (1, 1);
    \coordinate (B1) at (0, 2);
    \coordinate (B2) at (1, 2);
    \coordinate (B3) at (0, 3);
    \coordinate (B4) at (1, 3);
   \coordinate (C1) at (0, 4);
    \coordinate (C2) at (1, 4);
   \coordinate(D1) at (-1,2);
   \coordinate(D2) at (-1,3);
  \coordinate(D3) at (2,2);
 \coordinate(D4) at (2,3);

   \draw[very thick] (A1) -- (A2) ;
    \draw[very thick] (A1) -- (A3) ;
    \draw[very thick] (A2) -- (A4) ;
    \draw[very thick] (A3) -- (A4);
   \draw[very thick] (A3) -- (B1);
   \draw[very thick] (A4) -- (B2);
   \draw[very thick] (B1) -- (B2);
   \draw[very thick] (B3) -- (C1);
   \draw[very thick] (B4) -- (C2);
   \draw[very thick] (B3) -- (B1);
   \draw[very thick] (B4) -- (B2);
   \draw[very thick] (B3) -- (B4);
 \draw[very thick] (C1) -- (C2);    
  \draw[very thick] (D1) -- (D2); 
 \draw[very thick] (D3) -- (D4);
  \draw[very thick] (B1) -- (D1);
  \draw[very thick] (B2) -- (D3);
  \draw[very thick] (B3) -- (D2);
  \draw[very thick] (B4) -- (D4);

\node at (.5,.5) {3};
\node at (.5,1.5) {5};
\node at (.5,2.5) {5};
\node at (.5,3.5) {11};
\node at (-.5,2.5) {2};
\node at (1.5,2.5) {7};
\end{tikzpicture}

\begin{tikzpicture}[thick,scale=.5, ]
   \coordinate (A1) at (0,0);
   \coordinate (A2) at (1, 0);
    \coordinate (A3) at (0, 1);
    \coordinate (A4) at (1, 1);
    \coordinate (B1) at (0, 2);
    \coordinate (B2) at (1, 2);
    \coordinate (B3) at (0, 3);
    \coordinate (B4) at (1, 3);
   \coordinate (C1) at (0, 4);
    \coordinate (C2) at (1, 4);
   \coordinate(D1) at (-1,2);
   \coordinate(D2) at (-1,3);
  \coordinate(D3) at (2,2);
 \coordinate(D4) at (2,3);

   \draw[very thick] (A1) -- (A2);
    \draw[very thick] (A1) -- (A3);
    \draw[very thick] (A2) -- (A4);
    \draw[very thick] (A3) -- (A4);
   \draw[very thick] (A3) -- (B1);
   \draw[very thick] (A4) -- (B2);
   \draw[very thick] (B1) -- (B2);
   \draw[very thick] (B3) -- (C1);
   \draw[very thick] (B4) -- (C2);
   \draw[very thick] (B3) -- (B1);
   \draw[very thick] (B4) -- (B2);
   \draw[very thick] (B3) -- (B4);
 \draw[very thick] (C1) -- (C2);    
  \draw[very thick] (D1) -- (D2); 
 \draw[very thick] (D3) -- (D4);
  \draw[very thick] (B1) -- (D1);
  \draw[very thick] (B2) -- (D3);
  \draw[very thick] (B3) -- (D2);
  \draw[very thick] (B4) -- (D4);

\node at (.5,.5) {3};
\node at (.5,1.5) {7};
\node at (.5,2.5) {7};
\node at (.5,3.5) {11};
\node at (-.5,2.5) {5};
\node at (1.5,2.5) {2};
\end{tikzpicture}

\begin{tikzpicture}[thick,scale=.5, ]
   \coordinate (A1) at (0,0);
   \coordinate (A2) at (1, 0);
    \coordinate (A3) at (0, 1);
    \coordinate (A4) at (1, 1);
    \coordinate (B1) at (0, 2);
    \coordinate (B2) at (1, 2);
    \coordinate (B3) at (0, 3);
    \coordinate (B4) at (1, 3);
   \coordinate (C1) at (0, 4);
    \coordinate (C2) at (1, 4);
   \coordinate(D1) at (-1,2);
   \coordinate(D2) at (-1,3);
  \coordinate(D3) at (2,2);
 \coordinate(D4) at (2,3);

   \draw[very thick] (A1) -- (A2);
    \draw[very thick] (A1) -- (A3);
    \draw[very thick] (A2) -- (A4);
    \draw[very thick] (A3) -- (A4);
   \draw[very thick] (A3) -- (B1);
   \draw[very thick] (A4) -- (B2);
   \draw[very thick] (B1) -- (B2);
   \draw[very thick] (B3) -- (C1);
   \draw[very thick] (B4) -- (C2);
   \draw[very thick] (B3) -- (B1);
   \draw[very thick] (B4) -- (B2);
   \draw[very thick] (B3) -- (B4);
 \draw[very thick] (C1) -- (C2);    
  \draw[very thick] (D1) -- (D2); 
 \draw[very thick] (D3) -- (D4);
  \draw[very thick] (B1) -- (D1);
  \draw[very thick] (B2) -- (D3);
  \draw[very thick] (B3) -- (D2);
  \draw[very thick] (B4) -- (D4);

\node at (.5,.5) {3};
\node at (.5,1.5) {7};
\node at (.5,2.5) {11};
\node at (.5,3.5) {11};
\node at (-.5,2.5) {2};
\node at (1.5,2.5) {5};
\end{tikzpicture}

\end{multicols}
\end{center}

For $n=6$ there exists an Insanity puzzle with $18$ solutions and there are no Insanity puzzles with more than $18$ solutions. It is not true that for all $0 \leq m \leq 18$ there exists an Insanity puzzle with $m$ solutions. For $m= 5, 8, 10, 12, 14, 15, 16$, and $m=17$ there is no Insanity puzzle with 5 cubes and $m$ distinct solutions.

The next Insanity puzzle with six cubes has only one solution.

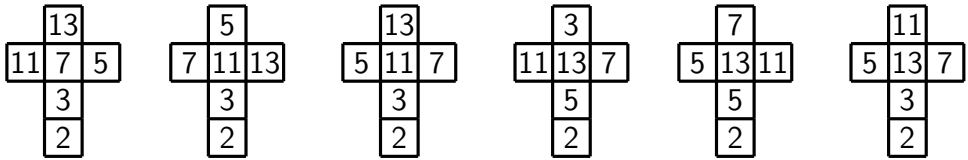
\begin{center}
 \begin{multicols}{6}

\begin{tikzpicture}[thick,scale=.5, ]
   \coordinate (A1) at (0,0);
   \coordinate (A2) at (1, 0);
    \coordinate (A3) at (0, 1);
    \coordinate (A4) at (1, 1);
    \coordinate (B1) at (0, 2);
    \coordinate (B2) at (1, 2);
    \coordinate (B3) at (0, 3);
    \coordinate (B4) at (1, 3);
   \coordinate (C1) at (0, 4);
    \coordinate (C2) at (1, 4);
   \coordinate(D1) at (-1,2);
   \coordinate(D2) at (-1,3);
  \coordinate(D3) at (2,2);
 \coordinate(D4) at (2,3);

   \draw[very thick] (A1) -- (A2);
    \draw[very thick] (A1) -- (A3);
    \draw[very thick] (A2) -- (A4);
    \draw[very thick] (A3) -- (A4);
   \draw[very thick] (A3) -- (B1);
   \draw[very thick] (A4) -- (B2);
   \draw[very thick] (B1) -- (B2);
   \draw[very thick] (B3) -- (C1);
   \draw[very thick] (B4) -- (C2);
   \draw[very thick] (B3) -- (B1);
   \draw[very thick] (B4) -- (B2);
   \draw[very thick] (B3) -- (B4);
 \draw[very thick] (C1) -- (C2);    
  \draw[very thick] (D1) -- (D2); 
 \draw[very thick] (D3) -- (D4);
  \draw[very thick] (B1) -- (D1);
  \draw[very thick] (B2) -- (D3);
  \draw[very thick] (B3) -- (D2);
  \draw[very thick] (B4) -- (D4);

\node at (.5,.5) {2};
\node at (.5,1.5) {3};
\node at (.5,2.5) {7};
\node at (.5,3.5) {13};
\node at (-.5,2.5) {11};
\node at (1.5,2.5) {5};
\end{tikzpicture}

\begin{tikzpicture}[thick,scale=.5,]
   \coordinate (A1) at (0,0);
   \coordinate (A2) at (1, 0);
    \coordinate (A3) at (0, 1);
    \coordinate (A4) at (1, 1);
    \coordinate (B1) at (0, 2);
    \coordinate (B2) at (1, 2);
    \coordinate (B3) at (0, 3);
    \coordinate (B4) at (1, 3);
   \coordinate (C1) at (0, 4);
    \coordinate (C2) at (1, 4);
   \coordinate(D1) at (-1,2);
   \coordinate(D2) at (-1,3);
  \coordinate(D3) at (2,2);
 \coordinate(D4) at (2,3);

   \draw[very thick] (A1) -- (A2);
    \draw[very thick] (A1) -- (A3);
    \draw[very thick] (A2) -- (A4);
    \draw[very thick] (A3) -- (A4);
   \draw[very thick] (A3) -- (B1);
   \draw[very thick] (A4) -- (B2);
   \draw[very thick] (B1) -- (B2);
   \draw[very thick] (B3) -- (C1);
   \draw[very thick] (B4) -- (C2);
   \draw[very thick] (B3) -- (B1);
   \draw[very thick] (B4) -- (B2);
   \draw[very thick] (B3) -- (B4);
 \draw[very thick] (C1) -- (C2);    
  \draw[very thick] (D1) -- (D2); 
 \draw[very thick] (D3) -- (D4);
  \draw[very thick] (B1) -- (D1);
  \draw[very thick] (B2) -- (D3);
  \draw[very thick] (B3) -- (D2);
  \draw[very thick] (B4) -- (D4);

\node at (.5,.5) {2};
\node at (.5,1.5) {3};
\node at (.5,2.5) {11};
\node at (.5,3.5) {5};
\node at (-.5,2.5) {7};
\node at (1.5,2.5) {13};

\end{tikzpicture}

\begin{tikzpicture}[thick,scale=.5,]
   \coordinate (A1) at (0,0);
   \coordinate (A2) at (1, 0);
    \coordinate (A3) at (0, 1);
    \coordinate (A4) at (1, 1);
    \coordinate (B1) at (0, 2);
    \coordinate (B2) at (1, 2);
    \coordinate (B3) at (0, 3);
    \coordinate (B4) at (1, 3);
   \coordinate (C1) at (0, 4);
    \coordinate (C2) at (1, 4);
   \coordinate(D1) at (-1,2);
   \coordinate(D2) at (-1,3);
  \coordinate(D3) at (2,2);
 \coordinate(D4) at (2,3);

   \draw[very thick] (A1) -- (A2) ;
    \draw[very thick] (A1) -- (A3) ;
    \draw[very thick] (A2) -- (A4) ;
    \draw[very thick] (A3) -- (A4);
   \draw[very thick] (A3) -- (B1);
   \draw[very thick] (A4) -- (B2);
   \draw[very thick] (B1) -- (B2);
   \draw[very thick] (B3) -- (C1);
   \draw[very thick] (B4) -- (C2);
   \draw[very thick] (B3) -- (B1);
   \draw[very thick] (B4) -- (B2);
   \draw[very thick] (B3) -- (B4);
 \draw[very thick] (C1) -- (C2);    
  \draw[very thick] (D1) -- (D2); 
 \draw[very thick] (D3) -- (D4);
  \draw[very thick] (B1) -- (D1);
  \draw[very thick] (B2) -- (D3);
  \draw[very thick] (B3) -- (D2);
  \draw[very thick] (B4) -- (D4);

\node at (.5,.5) {2};
\node at (.5,1.5) {3};
\node at (.5,2.5) {11};
\node at (.5,3.5) {13};
\node at (-.5,2.5) {5};
\node at (1.5,2.5) {7};
\end{tikzpicture}

\begin{tikzpicture}[thick,scale=.5, ]
   \coordinate (A1) at (0,0);
   \coordinate (A2) at (1, 0);
    \coordinate (A3) at (0, 1);
    \coordinate (A4) at (1, 1);
    \coordinate (B1) at (0, 2);
    \coordinate (B2) at (1, 2);
    \coordinate (B3) at (0, 3);
    \coordinate (B4) at (1, 3);
   \coordinate (C1) at (0, 4);
    \coordinate (C2) at (1, 4);
   \coordinate(D1) at (-1,2);
   \coordinate(D2) at (-1,3);
  \coordinate(D3) at (2,2);
 \coordinate(D4) at (2,3);

   \draw[very thick] (A1) -- (A2);
    \draw[very thick] (A1) -- (A3);
    \draw[very thick] (A2) -- (A4);
    \draw[very thick] (A3) -- (A4);
   \draw[very thick] (A3) -- (B1);
   \draw[very thick] (A4) -- (B2);
   \draw[very thick] (B1) -- (B2);
   \draw[very thick] (B3) -- (C1);
   \draw[very thick] (B4) -- (C2);
   \draw[very thick] (B3) -- (B1);
   \draw[very thick] (B4) -- (B2);
   \draw[very thick] (B3) -- (B4);
 \draw[very thick] (C1) -- (C2);    
  \draw[very thick] (D1) -- (D2); 
 \draw[very thick] (D3) -- (D4);
  \draw[very thick] (B1) -- (D1);
  \draw[very thick] (B2) -- (D3);
  \draw[very thick] (B3) -- (D2);
  \draw[very thick] (B4) -- (D4);

\node at (.5,.5) {2};
\node at (.5,1.5) {5};
\node at (.5,2.5) {13};
\node at (.5,3.5) {3};
\node at (-.5,2.5) {11};
\node at (1.5,2.5) {7};
\end{tikzpicture}

\begin{tikzpicture}[thick,scale=.5,]
   \coordinate (A1) at (0,0);
   \coordinate (A2) at (1, 0);
    \coordinate (A3) at (0, 1);
    \coordinate (A4) at (1, 1);
    \coordinate (B1) at (0, 2);
    \coordinate (B2) at (1, 2);
    \coordinate (B3) at (0, 3);
    \coordinate (B4) at (1, 3);
   \coordinate (C1) at (0, 4);
    \coordinate (C2) at (1, 4);
   \coordinate(D1) at (-1,2);
   \coordinate(D2) at (-1,3);
  \coordinate(D3) at (2,2);
 \coordinate(D4) at (2,3);

   \draw[very thick] (A1) -- (A2) ;
    \draw[very thick] (A1) -- (A3) ;
    \draw[very thick] (A2) -- (A4) ;
    \draw[very thick] (A3) -- (A4);
   \draw[very thick] (A3) -- (B1);
   \draw[very thick] (A4) -- (B2);
   \draw[very thick] (B1) -- (B2);
   \draw[very thick] (B3) -- (C1);
   \draw[very thick] (B4) -- (C2);
   \draw[very thick] (B3) -- (B1);
   \draw[very thick] (B4) -- (B2);
   \draw[very thick] (B3) -- (B4);
 \draw[very thick] (C1) -- (C2);    
  \draw[very thick] (D1) -- (D2); 
 \draw[very thick] (D3) -- (D4);
  \draw[very thick] (B1) -- (D1);
  \draw[very thick] (B2) -- (D3);
  \draw[very thick] (B3) -- (D2);
  \draw[very thick] (B4) -- (D4);

\node at (.5,.5) {2};
\node at (.5,1.5) {5};
\node at (.5,2.5) {13};
\node at (.5,3.5) {7};
\node at (-.5,2.5) {5};
\node at (1.5,2.5) {11};
\end{tikzpicture}

\begin{tikzpicture}[thick,scale=.5, ]
   \coordinate (A1) at (0,0);
   \coordinate (A2) at (1, 0);
    \coordinate (A3) at (0, 1);
    \coordinate (A4) at (1, 1);
    \coordinate (B1) at (0, 2);
    \coordinate (B2) at (1, 2);
    \coordinate (B3) at (0, 3);
    \coordinate (B4) at (1, 3);
   \coordinate (C1) at (0, 4);
    \coordinate (C2) at (1, 4);
   \coordinate(D1) at (-1,2);
   \coordinate(D2) at (-1,3);
  \coordinate(D3) at (2,2);
 \coordinate(D4) at (2,3);

   \draw[very thick] (A1) -- (A2);
    \draw[very thick] (A1) -- (A3);
    \draw[very thick] (A2) -- (A4);
    \draw[very thick] (A3) -- (A4);
   \draw[very thick] (A3) -- (B1);
   \draw[very thick] (A4) -- (B2);
   \draw[very thick] (B1) -- (B2);
   \draw[very thick] (B3) -- (C1);
   \draw[very thick] (B4) -- (C2);
   \draw[very thick] (B3) -- (B1);
   \draw[very thick] (B4) -- (B2);
   \draw[very thick] (B3) -- (B4);
 \draw[very thick] (C1) -- (C2);    
  \draw[very thick] (D1) -- (D2); 
 \draw[very thick] (D3) -- (D4);
  \draw[very thick] (B1) -- (D1);
  \draw[very thick] (B2) -- (D3);
  \draw[very thick] (B3) -- (D2);
  \draw[very thick] (B4) -- (D4);

\node at (.5,.5) {2};
\node at (.5,1.5) {3};
\node at (.5,2.5) {13};
\node at (.5,3.5) {11};
\node at (-.5,2.5) {5};
\node at (1.5,2.5) {7};
\end{tikzpicture}
\end{multicols}
\end{center}
 
 The next Insanity puzzle with five cubes has $18$ solutions.
 
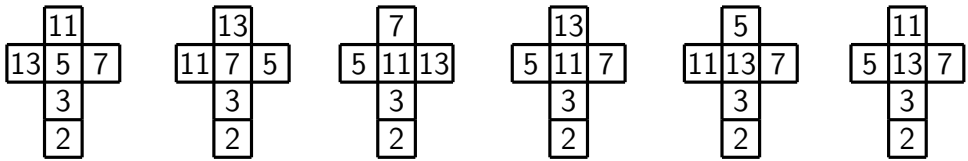
\begin{center}
 \begin{multicols}{6}

\begin{tikzpicture}[thick,scale=.5, ]
   \coordinate (A1) at (0,0);
   \coordinate (A2) at (1, 0);
    \coordinate (A3) at (0, 1);
    \coordinate (A4) at (1, 1);
    \coordinate (B1) at (0, 2);
    \coordinate (B2) at (1, 2);
    \coordinate (B3) at (0, 3);
    \coordinate (B4) at (1, 3);
   \coordinate (C1) at (0, 4);
    \coordinate (C2) at (1, 4);
   \coordinate(D1) at (-1,2);
   \coordinate(D2) at (-1,3);
  \coordinate(D3) at (2,2);
 \coordinate(D4) at (2,3);

   \draw[very thick] (A1) -- (A2);
    \draw[very thick] (A1) -- (A3);
    \draw[very thick] (A2) -- (A4);
    \draw[very thick] (A3) -- (A4);
   \draw[very thick] (A3) -- (B1);
   \draw[very thick] (A4) -- (B2);
   \draw[very thick] (B1) -- (B2);
   \draw[very thick] (B3) -- (C1);
   \draw[very thick] (B4) -- (C2);
   \draw[very thick] (B3) -- (B1);
   \draw[very thick] (B4) -- (B2);
   \draw[very thick] (B3) -- (B4);
 \draw[very thick] (C1) -- (C2);    
  \draw[very thick] (D1) -- (D2); 
 \draw[very thick] (D3) -- (D4);
  \draw[very thick] (B1) -- (D1);
  \draw[very thick] (B2) -- (D3);
  \draw[very thick] (B3) -- (D2);
  \draw[very thick] (B4) -- (D4);

\node at (.5,.5) {2};
\node at (.5,1.5) {3};
\node at (.5,2.5) {5};
\node at (.5,3.5) {11};
\node at (-.5,2.5) {13};
\node at (1.5,2.5) {7};
\end{tikzpicture}

\begin{tikzpicture}[thick,scale=.5,]
   \coordinate (A1) at (0,0);
   \coordinate (A2) at (1, 0);
    \coordinate (A3) at (0, 1);
    \coordinate (A4) at (1, 1);
    \coordinate (B1) at (0, 2);
    \coordinate (B2) at (1, 2);
    \coordinate (B3) at (0, 3);
    \coordinate (B4) at (1, 3);
   \coordinate (C1) at (0, 4);
    \coordinate (C2) at (1, 4);
   \coordinate(D1) at (-1,2);
   \coordinate(D2) at (-1,3);
  \coordinate(D3) at (2,2);
 \coordinate(D4) at (2,3);

   \draw[very thick] (A1) -- (A2);
    \draw[very thick] (A1) -- (A3);
    \draw[very thick] (A2) -- (A4);
    \draw[very thick] (A3) -- (A4);
   \draw[very thick] (A3) -- (B1);
   \draw[very thick] (A4) -- (B2);
   \draw[very thick] (B1) -- (B2);
   \draw[very thick] (B3) -- (C1);
   \draw[very thick] (B4) -- (C2);
   \draw[very thick] (B3) -- (B1);
   \draw[very thick] (B4) -- (B2);
   \draw[very thick] (B3) -- (B4);
 \draw[very thick] (C1) -- (C2);    
  \draw[very thick] (D1) -- (D2); 
 \draw[very thick] (D3) -- (D4);
  \draw[very thick] (B1) -- (D1);
  \draw[very thick] (B2) -- (D3);
  \draw[very thick] (B3) -- (D2);
  \draw[very thick] (B4) -- (D4);

\node at (.5,.5) {2};
\node at (.5,1.5) {3};
\node at (.5,2.5) {7};
\node at (.5,3.5) {13};
\node at (-.5,2.5) {11};
\node at (1.5,2.5) {5};

\end{tikzpicture}

\begin{tikzpicture}[thick,scale=.5,]
   \coordinate (A1) at (0,0);
   \coordinate (A2) at (1, 0);
    \coordinate (A3) at (0, 1);
    \coordinate (A4) at (1, 1);
    \coordinate (B1) at (0, 2);
    \coordinate (B2) at (1, 2);
    \coordinate (B3) at (0, 3);
    \coordinate (B4) at (1, 3);
   \coordinate (C1) at (0, 4);
    \coordinate (C2) at (1, 4);
   \coordinate(D1) at (-1,2);
   \coordinate(D2) at (-1,3);
  \coordinate(D3) at (2,2);
 \coordinate(D4) at (2,3);

   \draw[very thick] (A1) -- (A2) ;
    \draw[very thick] (A1) -- (A3) ;
    \draw[very thick] (A2) -- (A4) ;
    \draw[very thick] (A3) -- (A4);
   \draw[very thick] (A3) -- (B1);
   \draw[very thick] (A4) -- (B2);
   \draw[very thick] (B1) -- (B2);
   \draw[very thick] (B3) -- (C1);
   \draw[very thick] (B4) -- (C2);
   \draw[very thick] (B3) -- (B1);
   \draw[very thick] (B4) -- (B2);
   \draw[very thick] (B3) -- (B4);
 \draw[very thick] (C1) -- (C2);    
  \draw[very thick] (D1) -- (D2); 
 \draw[very thick] (D3) -- (D4);
  \draw[very thick] (B1) -- (D1);
  \draw[very thick] (B2) -- (D3);
  \draw[very thick] (B3) -- (D2);
  \draw[very thick] (B4) -- (D4);

\node at (.5,.5) {2};
\node at (.5,1.5) {3};
\node at (.5,2.5) {11};
\node at (.5,3.5) {7};
\node at (-.5,2.5) {5};
\node at (1.5,2.5) {13};
\end{tikzpicture}

\begin{tikzpicture}[thick,scale=.5, ]
   \coordinate (A1) at (0,0);
   \coordinate (A2) at (1, 0);
    \coordinate (A3) at (0, 1);
    \coordinate (A4) at (1, 1);
    \coordinate (B1) at (0, 2);
    \coordinate (B2) at (1, 2);
    \coordinate (B3) at (0, 3);
    \coordinate (B4) at (1, 3);
   \coordinate (C1) at (0, 4);
    \coordinate (C2) at (1, 4);
   \coordinate(D1) at (-1,2);
   \coordinate(D2) at (-1,3);
  \coordinate(D3) at (2,2);
 \coordinate(D4) at (2,3);

   \draw[very thick] (A1) -- (A2);
    \draw[very thick] (A1) -- (A3);
    \draw[very thick] (A2) -- (A4);
    \draw[very thick] (A3) -- (A4);
   \draw[very thick] (A3) -- (B1);
   \draw[very thick] (A4) -- (B2);
   \draw[very thick] (B1) -- (B2);
   \draw[very thick] (B3) -- (C1);
   \draw[very thick] (B4) -- (C2);
   \draw[very thick] (B3) -- (B1);
   \draw[very thick] (B4) -- (B2);
   \draw[very thick] (B3) -- (B4);
 \draw[very thick] (C1) -- (C2);    
  \draw[very thick] (D1) -- (D2); 
 \draw[very thick] (D3) -- (D4);
  \draw[very thick] (B1) -- (D1);
  \draw[very thick] (B2) -- (D3);
  \draw[very thick] (B3) -- (D2);
  \draw[very thick] (B4) -- (D4);

\node at (.5,.5) {2};
\node at (.5,1.5) {3};
\node at (.5,2.5) {11};
\node at (.5,3.5) {13};
\node at (-.5,2.5) {5};
\node at (1.5,2.5) {7};
\end{tikzpicture}

\begin{tikzpicture}[thick,scale=.5,]
   \coordinate (A1) at (0,0);
   \coordinate (A2) at (1, 0);
    \coordinate (A3) at (0, 1);
    \coordinate (A4) at (1, 1);
    \coordinate (B1) at (0, 2);
    \coordinate (B2) at (1, 2);
    \coordinate (B3) at (0, 3);
    \coordinate (B4) at (1, 3);
   \coordinate (C1) at (0, 4);
    \coordinate (C2) at (1, 4);
   \coordinate(D1) at (-1,2);
   \coordinate(D2) at (-1,3);
  \coordinate(D3) at (2,2);
 \coordinate(D4) at (2,3);

   \draw[very thick] (A1) -- (A2) ;
    \draw[very thick] (A1) -- (A3) ;
    \draw[very thick] (A2) -- (A4) ;
    \draw[very thick] (A3) -- (A4);
   \draw[very thick] (A3) -- (B1);
   \draw[very thick] (A4) -- (B2);
   \draw[very thick] (B1) -- (B2);
   \draw[very thick] (B3) -- (C1);
   \draw[very thick] (B4) -- (C2);
   \draw[very thick] (B3) -- (B1);
   \draw[very thick] (B4) -- (B2);
   \draw[very thick] (B3) -- (B4);
 \draw[very thick] (C1) -- (C2);    
  \draw[very thick] (D1) -- (D2); 
 \draw[very thick] (D3) -- (D4);
  \draw[very thick] (B1) -- (D1);
  \draw[very thick] (B2) -- (D3);
  \draw[very thick] (B3) -- (D2);
  \draw[very thick] (B4) -- (D4);

\node at (.5,.5) {2};
\node at (.5,1.5) {3};
\node at (.5,2.5) {13};
\node at (.5,3.5) {5};
\node at (-.5,2.5) {11};
\node at (1.5,2.5) {7};
\end{tikzpicture}

\begin{tikzpicture}[thick,scale=.5, ]
   \coordinate (A1) at (0,0);
   \coordinate (A2) at (1, 0);
    \coordinate (A3) at (0, 1);
    \coordinate (A4) at (1, 1);
    \coordinate (B1) at (0, 2);
    \coordinate (B2) at (1, 2);
    \coordinate (B3) at (0, 3);
    \coordinate (B4) at (1, 3);
   \coordinate (C1) at (0, 4);
    \coordinate (C2) at (1, 4);
   \coordinate(D1) at (-1,2);
   \coordinate(D2) at (-1,3);
  \coordinate(D3) at (2,2);
 \coordinate(D4) at (2,3);

   \draw[very thick] (A1) -- (A2);
    \draw[very thick] (A1) -- (A3);
    \draw[very thick] (A2) -- (A4);
    \draw[very thick] (A3) -- (A4);
   \draw[very thick] (A3) -- (B1);
   \draw[very thick] (A4) -- (B2);
   \draw[very thick] (B1) -- (B2);
   \draw[very thick] (B3) -- (C1);
   \draw[very thick] (B4) -- (C2);
   \draw[very thick] (B3) -- (B1);
   \draw[very thick] (B4) -- (B2);
   \draw[very thick] (B3) -- (B4);
 \draw[very thick] (C1) -- (C2);    
  \draw[very thick] (D1) -- (D2); 
 \draw[very thick] (D3) -- (D4);
  \draw[very thick] (B1) -- (D1);
  \draw[very thick] (B2) -- (D3);
  \draw[very thick] (B3) -- (D2);
  \draw[very thick] (B4) -- (D4);

\node at (.5,.5) {2};
\node at (.5,1.5) {3};
\node at (.5,2.5) {13};
\node at (.5,3.5) {11};
\node at (-.5,2.5) {5};
\node at (1.5,2.5) {7};
\end{tikzpicture}
\end{multicols}
\end{center}

\section*{The Mutando of Insanity}

The Mutando is a puzzle designed  by E. Kunzell in 1997. It was commercialise in 2000 by Ingo Uhl in 2000. This puzzle consists of four cubes coloured with four colours with a structure analogous to the next four cubes.

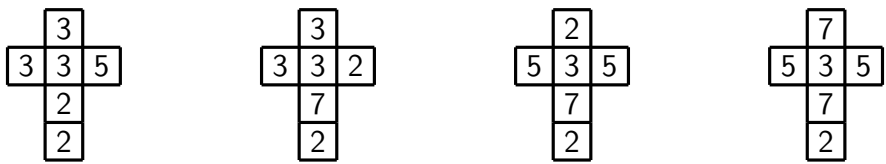
\begin{center}
 \begin{multicols}{4}
\begin{tikzpicture}[thick,scale=.5, ]
   \coordinate (A1) at (0,0);
   \coordinate (A2) at (1, 0);
    \coordinate (A3) at (0, 1);
    \coordinate (A4) at (1, 1);
    \coordinate (B1) at (0, 2);
    \coordinate (B2) at (1, 2);
    \coordinate (B3) at (0, 3);
    \coordinate (B4) at (1, 3);
   \coordinate (C1) at (0, 4);
    \coordinate (C2) at (1, 4);
   \coordinate(D1) at (-1,2);
   \coordinate(D2) at (-1,3);
  \coordinate(D3) at (2,2);
 \coordinate(D4) at (2,3);

   \draw[very thick] (A1) -- (A2);
    \draw[very thick] (A1) -- (A3);
    \draw[very thick] (A2) -- (A4);
    \draw[very thick] (A3) -- (A4);
   \draw[very thick] (A3) -- (B1);
   \draw[very thick] (A4) -- (B2);
   \draw[very thick] (B1) -- (B2);
   \draw[very thick] (B3) -- (C1);
   \draw[very thick] (B4) -- (C2);
   \draw[very thick] (B3) -- (B1);
   \draw[very thick] (B4) -- (B2);
   \draw[very thick] (B3) -- (B4);
 \draw[very thick] (C1) -- (C2);    
  \draw[very thick] (D1) -- (D2); 
 \draw[very thick] (D3) -- (D4);
  \draw[very thick] (B1) -- (D1);
  \draw[very thick] (B2) -- (D3);
  \draw[very thick] (B3) -- (D2);
  \draw[very thick] (B4) -- (D4);

\node at (.5,.5) {2};
\node at (.5,1.5) {2};
\node at (.5,2.5) {3};
\node at (.5,3.5) {3};
\node at (-.5,2.5) {3};
\node at (1.5,2.5) {5};
\end{tikzpicture}

\begin{tikzpicture}[thick,scale=.5,]
   \coordinate (A1) at (0,0);
   \coordinate (A2) at (1, 0);
    \coordinate (A3) at (0, 1);
    \coordinate (A4) at (1, 1);
    \coordinate (B1) at (0, 2);
    \coordinate (B2) at (1, 2);
    \coordinate (B3) at (0, 3);
    \coordinate (B4) at (1, 3);
   \coordinate (C1) at (0, 4);
    \coordinate (C2) at (1, 4);
   \coordinate(D1) at (-1,2);
   \coordinate(D2) at (-1,3);
  \coordinate(D3) at (2,2);
 \coordinate(D4) at (2,3);

   \draw[very thick] (A1) -- (A2);
    \draw[very thick] (A1) -- (A3);
    \draw[very thick] (A2) -- (A4);
    \draw[very thick] (A3) -- (A4);
   \draw[very thick] (A3) -- (B1);
   \draw[very thick] (A4) -- (B2);
   \draw[very thick] (B1) -- (B2);
   \draw[very thick] (B3) -- (C1);
   \draw[very thick] (B4) -- (C2);
   \draw[very thick] (B3) -- (B1);
   \draw[very thick] (B4) -- (B2);
   \draw[very thick] (B3) -- (B4);
 \draw[very thick] (C1) -- (C2);    
  \draw[very thick] (D1) -- (D2); 
 \draw[very thick] (D3) -- (D4);
  \draw[very thick] (B1) -- (D1);
  \draw[very thick] (B2) -- (D3);
  \draw[very thick] (B3) -- (D2);
  \draw[very thick] (B4) -- (D4);

\node at (.5,.5) {2};
\node at (.5,1.5) {7};
\node at (.5,2.5) {3};
\node at (.5,3.5) {3};
\node at (-.5,2.5) {3};
\node at (1.5,2.5) {2};

\end{tikzpicture}

\begin{tikzpicture}[thick,scale=.5,]
   \coordinate (A1) at (0,0);
   \coordinate (A2) at (1, 0);
    \coordinate (A3) at (0, 1);
    \coordinate (A4) at (1, 1);
    \coordinate (B1) at (0, 2);
    \coordinate (B2) at (1, 2);
    \coordinate (B3) at (0, 3);
    \coordinate (B4) at (1, 3);
   \coordinate (C1) at (0, 4);
    \coordinate (C2) at (1, 4);
   \coordinate(D1) at (-1,2);
   \coordinate(D2) at (-1,3);
  \coordinate(D3) at (2,2);
 \coordinate(D4) at (2,3);

   \draw[very thick] (A1) -- (A2) ;
    \draw[very thick] (A1) -- (A3) ;
    \draw[very thick] (A2) -- (A4) ;
    \draw[very thick] (A3) -- (A4);
   \draw[very thick] (A3) -- (B1);
   \draw[very thick] (A4) -- (B2);
   \draw[very thick] (B1) -- (B2);
   \draw[very thick] (B3) -- (C1);
   \draw[very thick] (B4) -- (C2);
   \draw[very thick] (B3) -- (B1);
   \draw[very thick] (B4) -- (B2);
   \draw[very thick] (B3) -- (B4);
 \draw[very thick] (C1) -- (C2);    
  \draw[very thick] (D1) -- (D2); 
 \draw[very thick] (D3) -- (D4);
  \draw[very thick] (B1) -- (D1);
  \draw[very thick] (B2) -- (D3);
  \draw[very thick] (B3) -- (D2);
  \draw[very thick] (B4) -- (D4);

\node at (.5,.5) {2};
\node at (.5,1.5) {7};
\node at (.5,2.5) {3};
\node at (.5,3.5) {2};
\node at (-.5,2.5) {5};
\node at (1.5,2.5) {5};
\end{tikzpicture}

\begin{tikzpicture}[thick,scale=.5, ]
   \coordinate (A1) at (0,0);
   \coordinate (A2) at (1, 0);
    \coordinate (A3) at (0, 1);
    \coordinate (A4) at (1, 1);
    \coordinate (B1) at (0, 2);
    \coordinate (B2) at (1, 2);
    \coordinate (B3) at (0, 3);
    \coordinate (B4) at (1, 3);
   \coordinate (C1) at (0, 4);
    \coordinate (C2) at (1, 4);
   \coordinate(D1) at (-1,2);
   \coordinate(D2) at (-1,3);
  \coordinate(D3) at (2,2);
 \coordinate(D4) at (2,3);

   \draw[very thick] (A1) -- (A2);
    \draw[very thick] (A1) -- (A3);
    \draw[very thick] (A2) -- (A4);
    \draw[very thick] (A3) -- (A4);
   \draw[very thick] (A3) -- (B1);
   \draw[very thick] (A4) -- (B2);
   \draw[very thick] (B1) -- (B2);
   \draw[very thick] (B3) -- (C1);
   \draw[very thick] (B4) -- (C2);
   \draw[very thick] (B3) -- (B1);
   \draw[very thick] (B4) -- (B2);
   \draw[very thick] (B3) -- (B4);
 \draw[very thick] (C1) -- (C2);    
  \draw[very thick] (D1) -- (D2); 
 \draw[very thick] (D3) -- (D4);
  \draw[very thick] (B1) -- (D1);
  \draw[very thick] (B2) -- (D3);
  \draw[very thick] (B3) -- (D2);
  \draw[very thick] (B4) -- (D4);

\node at (.5,.5) {2};
\node at (.5,1.5) {7};
\node at (.5,2.5) {3};
\node at (.5,3.5) {7};
\node at (-.5,2.5) {5};
\node at (1.5,2.5) {5};
\end{tikzpicture}
\end{multicols}
\end{center}
The Mutando is not an Insanity puzzle because it does not have the four labels present in the four cubes. The two puzzles that the Mutando asks to solve are:
\begin{itemize}
\item Puzzle 1: With the four cubes, form a $4 \times 1 \times 1$ prism in such a way that in each long face of the prism the four labels are present. 
\item Puzzle 2: With the four cubes form a $2 \times 2 \times 1$ prism in such a way that in each face of the prism has all square faces with the same label. 
\end{itemize}

Observe that Puzzle 1 is what we ask from an Insanity puzzle to be solved. It is a natural question to ask if there exists an Insanity puzzle (all labels present in each one of the cubes) for $n=4$ such that it is possible to solve Puzzle1 and Puzzle 2 as defined above.  

Based on the same mathematical model of prime product matrices we were able to find an Insanity puzzle that has a solution to  Puzzle 1 and Puzzle 2. We have already described in previous section the mathematical model that we have used to solve Puzzle 1. With the same mathematical structure, but representing with a product of two primes the labeled faces of the cube that share an edge (the dimensions of the matrices were not the same because there are 12 pair of faces of a cube sharing and edge instead of three pairs of opposite faces), we were also able to solve Puzzle 2.

This Insanity puzzle was presented in the G4G 12 as a Gift Exchange and the configuration of the four cubes is:

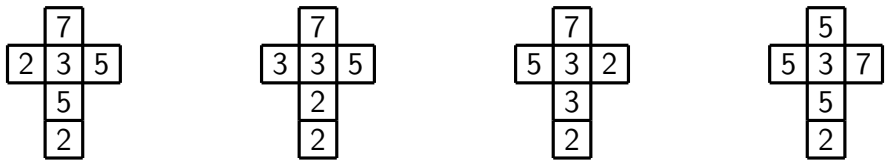
\begin{center}
 \begin{multicols}{4}
\begin{tikzpicture}[thick,scale=.5, ]
   \coordinate (A1) at (0,0);
   \coordinate (A2) at (1, 0);
    \coordinate (A3) at (0, 1);
    \coordinate (A4) at (1, 1);
    \coordinate (B1) at (0, 2);
    \coordinate (B2) at (1, 2);
    \coordinate (B3) at (0, 3);
    \coordinate (B4) at (1, 3);
   \coordinate (C1) at (0, 4);
    \coordinate (C2) at (1, 4);
   \coordinate(D1) at (-1,2);
   \coordinate(D2) at (-1,3);
  \coordinate(D3) at (2,2);
 \coordinate(D4) at (2,3);

   \draw[very thick] (A1) -- (A2);
    \draw[very thick] (A1) -- (A3);
    \draw[very thick] (A2) -- (A4);
    \draw[very thick] (A3) -- (A4);
   \draw[very thick] (A3) -- (B1);
   \draw[very thick] (A4) -- (B2);
   \draw[very thick] (B1) -- (B2);
   \draw[very thick] (B3) -- (C1);
   \draw[very thick] (B4) -- (C2);
   \draw[very thick] (B3) -- (B1);
   \draw[very thick] (B4) -- (B2);
   \draw[very thick] (B3) -- (B4);
 \draw[very thick] (C1) -- (C2);    
  \draw[very thick] (D1) -- (D2); 
 \draw[very thick] (D3) -- (D4);
  \draw[very thick] (B1) -- (D1);
  \draw[very thick] (B2) -- (D3);
  \draw[very thick] (B3) -- (D2);
  \draw[very thick] (B4) -- (D4);

\node at (.5,.5) {2};
\node at (.5,1.5) {5};
\node at (.5,2.5) {3};
\node at (.5,3.5) {7};
\node at (-.5,2.5) {2};
\node at (1.5,2.5) {5};
\end{tikzpicture}

\begin{tikzpicture}[thick,scale=.5,]
   \coordinate (A1) at (0,0);
   \coordinate (A2) at (1, 0);
    \coordinate (A3) at (0, 1);
    \coordinate (A4) at (1, 1);
    \coordinate (B1) at (0, 2);
    \coordinate (B2) at (1, 2);
    \coordinate (B3) at (0, 3);
    \coordinate (B4) at (1, 3);
   \coordinate (C1) at (0, 4);
    \coordinate (C2) at (1, 4);
   \coordinate(D1) at (-1,2);
   \coordinate(D2) at (-1,3);
  \coordinate(D3) at (2,2);
 \coordinate(D4) at (2,3);

   \draw[very thick] (A1) -- (A2);
    \draw[very thick] (A1) -- (A3);
    \draw[very thick] (A2) -- (A4);
    \draw[very thick] (A3) -- (A4);
   \draw[very thick] (A3) -- (B1);
   \draw[very thick] (A4) -- (B2);
   \draw[very thick] (B1) -- (B2);
   \draw[very thick] (B3) -- (C1);
   \draw[very thick] (B4) -- (C2);
   \draw[very thick] (B3) -- (B1);
   \draw[very thick] (B4) -- (B2);
   \draw[very thick] (B3) -- (B4);
 \draw[very thick] (C1) -- (C2);    
  \draw[very thick] (D1) -- (D2); 
 \draw[very thick] (D3) -- (D4);
  \draw[very thick] (B1) -- (D1);
  \draw[very thick] (B2) -- (D3);
  \draw[very thick] (B3) -- (D2);
  \draw[very thick] (B4) -- (D4);

\node at (.5,.5) {2};
\node at (.5,1.5) {2};
\node at (.5,2.5) {3};
\node at (.5,3.5) {7};
\node at (-.5,2.5) {3};
\node at (1.5,2.5) {5};

\end{tikzpicture}

\begin{tikzpicture}[thick,scale=.5,]
   \coordinate (A1) at (0,0);
   \coordinate (A2) at (1, 0);
    \coordinate (A3) at (0, 1);
    \coordinate (A4) at (1, 1);
    \coordinate (B1) at (0, 2);
    \coordinate (B2) at (1, 2);
    \coordinate (B3) at (0, 3);
    \coordinate (B4) at (1, 3);
   \coordinate (C1) at (0, 4);
    \coordinate (C2) at (1, 4);
   \coordinate(D1) at (-1,2);
   \coordinate(D2) at (-1,3);
  \coordinate(D3) at (2,2);
 \coordinate(D4) at (2,3);

   \draw[very thick] (A1) -- (A2) ;
    \draw[very thick] (A1) -- (A3) ;
    \draw[very thick] (A2) -- (A4) ;
    \draw[very thick] (A3) -- (A4);
   \draw[very thick] (A3) -- (B1);
   \draw[very thick] (A4) -- (B2);
   \draw[very thick] (B1) -- (B2);
   \draw[very thick] (B3) -- (C1);
   \draw[very thick] (B4) -- (C2);
   \draw[very thick] (B3) -- (B1);
   \draw[very thick] (B4) -- (B2);
   \draw[very thick] (B3) -- (B4);
 \draw[very thick] (C1) -- (C2);    
  \draw[very thick] (D1) -- (D2); 
 \draw[very thick] (D3) -- (D4);
  \draw[very thick] (B1) -- (D1);
  \draw[very thick] (B2) -- (D3);
  \draw[very thick] (B3) -- (D2);
  \draw[very thick] (B4) -- (D4);

\node at (.5,.5) {2};
\node at (.5,1.5) {3};
\node at (.5,2.5) {3};
\node at (.5,3.5) {7};
\node at (-.5,2.5) {5};
\node at (1.5,2.5) {2};
\end{tikzpicture}

\begin{tikzpicture}[thick,scale=.5, ]
   \coordinate (A1) at (0,0);
   \coordinate (A2) at (1, 0);
    \coordinate (A3) at (0, 1);
    \coordinate (A4) at (1, 1);
    \coordinate (B1) at (0, 2);
    \coordinate (B2) at (1, 2);
    \coordinate (B3) at (0, 3);
    \coordinate (B4) at (1, 3);
   \coordinate (C1) at (0, 4);
    \coordinate (C2) at (1, 4);
   \coordinate(D1) at (-1,2);
   \coordinate(D2) at (-1,3);
  \coordinate(D3) at (2,2);
 \coordinate(D4) at (2,3);

   \draw[very thick] (A1) -- (A2);
    \draw[very thick] (A1) -- (A3);
    \draw[very thick] (A2) -- (A4);
    \draw[very thick] (A3) -- (A4);
   \draw[very thick] (A3) -- (B1);
   \draw[very thick] (A4) -- (B2);
   \draw[very thick] (B1) -- (B2);
   \draw[very thick] (B3) -- (C1);
   \draw[very thick] (B4) -- (C2);
   \draw[very thick] (B3) -- (B1);
   \draw[very thick] (B4) -- (B2);
   \draw[very thick] (B3) -- (B4);
 \draw[very thick] (C1) -- (C2);    
  \draw[very thick] (D1) -- (D2); 
 \draw[very thick] (D3) -- (D4);
  \draw[very thick] (B1) -- (D1);
  \draw[very thick] (B2) -- (D3);
  \draw[very thick] (B3) -- (D2);
  \draw[very thick] (B4) -- (D4);

\node at (.5,.5) {2};
\node at (.5,1.5) {5};
\node at (.5,2.5) {3};
\node at (.5,3.5) {5};
\node at (-.5,2.5) {5};
\node at (1.5,2.5) {7};
\end{tikzpicture}
\end{multicols}
\end{center}


\section*{Acknowledgments }
I would like to thank Martin Gardner for the legacy that he had let us in very different ways. Part of his legacy is G4G. I was able to assist to the G4G 12 thanks to Colm Mulcahy invitation and the generous financial support that the G4G offers to students. At the event I had the opportunity to talked about this puzzles with Rik Van Grol and David Singmaster. I thank both for the interesting conversations and the information that they generously gave me about Insanity puzzles.  
\nocite{*}
\bibliographystyle{plain}
\bibliography{MIB}

\end{document}